\documentclass[]{gOMS2e}

\usepackage{psfrag}
\usepackage{rotating}
\usepackage[usenames]{color}
\usepackage{graphicx}
\graphicspath{{figures/}}
\usepackage{epsfig}    
\usepackage{amssymb}   
\usepackage[font=footnotesize,caption=false]{subfig}
\usepackage{fixltx2e}
\usepackage{array}
\usepackage{dblfloatfix}

\usepackage{url}
\usepackage{tikz}

\usepackage{algorithm}     
\usepackage{algpseudocode} 

\DeclareMathOperator*{\minimize}{minimize}

\DeclareMathOperator*{\minimum}{min}
\DeclareMathOperator*{\subject}{subject\ to}

\DeclareMathOperator*{\maximum}{max}
\DeclareMathOperator*{\argmin}{arg\ min}

\DeclareMathOperator*{\diag}{diag}
\DeclareMathOperator*{\relint}{rel\ int}
\DeclareMathOperator*{\ridom}{ri\ dom}
\DeclareMathOperator*{\dom}{dom}

\DeclareMathOperator*{\infimum}{inf}

\DeclareMathOperator*{\find}{find}
\DeclareMathOperator*{\adj}{adj}

\newcounter{remcount}
\newtheorem{rem}[remcount]{Remark}

\definecolor{red}{rgb}{1,0,0}
\newcommand{\etal}{\emph{et~al.}}

\newcommand{\dist}{\mathrm{dist}}

\begin{document}
\doi{10.1080/1055.6788.YYYY.xxxxxx}
\issn{1029-4937}
\issnp{1055-6788}
\jvol{00} \jnum{00} \jyear{2013} \jmonth{January}

\markboth{Taylor \& Francis and I.T. Consultant}{Optimization Methods and Software}

\title{Distributed Solutions for Loosely Coupled Feasibility Problems Using Proximal Splitting Methods$^*$\thanks{$^*$This work has been supported by the Swedish Department of Education within the ELLIIT project.}}


\author{Sina Khoshfetrat Pakazad$^{1}$\thanks{$^{1}$ Sina Khoshfetrat Pakazad and Anders Hansson are with the Division of Automatic Control, Department of Electrical Engineering, Link\"oping University, Sweden. Email: \{sina.kh.pa, hansson\}@isy.liu.se.}, Martin S. Andersen$^{2}$ \thanks{$^{2}$Martin S. Andersen is with the Department of Applied Mathematics and Computer Science, Technical University of Denmark. Email:  mskan@dtu.dk.} and Anders Hansson$^{1}$}

\maketitle


\begin{abstract}
In this paper, we consider convex feasibility problems where the underlying sets are loosely coupled, and we propose several algorithms to solve such problems in a distributed manner. These algorithms are obtained by applying proximal splitting methods to convex minimization reformulations of convex feasibility problems. We also put forth  distributed convergence tests which enable us to establish feasibility or infeasibility of the problem distributedly, and we provide convergence rate results. Under the assumption that the problem is feasible and boundedly linearly regular, these convergence results are given in terms of the distance of the iterates to the feasible set, which are similar to those of classical projection methods. In case the feasibility problem is infeasible, we provide convergence rate results that concern the convergence of certain error-bounds.

\begin{classcode}G.1.6; G.1.10; G.2.2; I.1.2 \end{classcode}\bigskip

\end{abstract}

\section{Introduction}\label{sec:Intro}
A convex feasibility problem (CFP), corresponds to the problem of finding an element in the intersection of, say, $N$ non-empty convex sets ($\mathcal C_i$), i.e., $ x \in \bigcap_{i=1}^{N} \mathcal C_i$. Such problems appear in many fields of engineering and science, e.g., image reconstruction, best approximation theory, analysis of networked systems \cite{kho:13,bau:96,her:95}, and have been studied thoroughly in the literature, e.g., see \cite{bau:96,bet:03,bau:93,bau:94,gub:99}. Many of the algorithms designed for solving CFPs, rely on projections onto the individual sets, and are referred to as projection methods. The behavior and convergence properties of such algorithms have been well-studied, \cite{bau:96,bet:03,bau:93,bau:94}, both when the CFP is feasible or  infeasible. For a thorough survey of such algorithms refer to \cite{bau:96}. In this paper, we focus on CFPs where each constraint defining a set $\mathcal C_i$ in the problem is only dependent on a subset of components of a variable, $x$. We also assume that the number of variables that appear jointly in the descriptions of every two constraint sets $\mathcal C_i$ and $\mathcal C_j$ ($i \neq j$) is small. We refer to such CFPs as loosely coupled, and we intend to develop distributed algorithms for solving these problems efficiently. Employing classical projection algorithms for solving such CFPs, while neglecting the underlying structure in the problem, is inefficient and has been shown to be extremely slow, \cite{kho:13}. In order to boost the performance of these algorithms,  the structure in the CFP needs to be exploited. This can be done by using ideas from \cite{boyd:11}, \cite{ber:97} and \cite{kiw:99}, and it results in a reformulation of the problem as an equivalent feasibility problem in product space. Similar formulations of CFPs are also proposed in \cite{dep:01,pie:84}. This product-space formulation can be solved using well-known projection methods, e.g., von Neumann's and Dykstra's alternating projections. Although such algorithms have been shown to be effective and their convergence properties are well-studied, to the best knowledge of the authors, rate of convergence of such methods are mainly available when the CFP is feasible, \cite{bet:03}. Hence, in order to provide convergence rate results for projection-based methods even when the problem is infeasible, we approach the problem in another way.

The product-space reformulation of loosely coupled CFPs can also be written as a convex minimization problem, \cite{bet:03,but:00,dep:85,cen:94}. Authors in \cite{bet:03}, define a convex minimization reformulation of CFPs and consider the use of gradient projection algorithms (GPA) for solving minimization problem. They then utilize the convergence properties of GPA to provide convergence results for the resulting algorithm. Furthermore, they show that in case the problem is infeasible, certain error bounds converge to non-zero values with $\mathcal O(1/\sqrt{k})$ rate of convergence. Prior to \cite{bet:03}, such approaches have also been used for similar purposes, using the so-called subgradient methods, \cite{bau:94,pol:69,dos:87,cen:82}. In this paper we also deal with convex minimization reformulations of CFPs. This in turn allows us to employ proximal splitting methods (first order methods) for solving minimization reformulations of CFPs, which result in several distributed projection-based algorithms for solving loosely coupled CFPs. Some of these algorithms are similar to already existing well-known projection methods. However, several of the proposed algorithms can be considered as generalizations of classical projection methods, \cite{cen:94}. Furthermore, if the minimization formulation of a CFP is well-defined even when it is infeasible, the convergence properties, and particularly the rate of convergence of the utilized proximal methods can be used to establish convergence results for the newly generated algorithms. This allows us to provide rate of convergence results for such algorithms even when the CFP is infeasible.

The contributions of this paper are as follows. In this paper, we
\begin{itemize}
\item propose several distributed algorithms for solving loosely coupled CFPs;
\item establish local convergence tests for these algorithms, which enable us to detect arrival at a feasible solution or infeasibility of the problem in a distributed manner with minimal communication;
\item also provide convergence rate results for the proposed algorithms for the case the CFP is either feasible or infeasible. For the case the problem is feasible, these results are given in terms of the distance of the iterates to the feasible set. This enables us to provide a unified treatment of the convergence rate analysis of these algorithms and the classic projection methods. In case of an infeasible problem, the result is a bound of the rate of convergence of a norm of a certain residual to a non-zero constant.
\end{itemize}
The performance of the proposed algorithms are compared using numerical examples.
\subsection*{Outline}
The paper is organized as follows. In Section \ref{sec:Decomposition},
we provide a formal description of loosely coupled CFPs  and describe different approaches for formulating and solving such problems. Particularly we discuss minimization reformulations of coupled CFPs in Section~\ref{sec:Convex}. A brief description of the most commonly used proximal splitting methods is given in Section~\ref{sec:prox}. These methods are then applied to convex minimization reformulations of the coupled CFP, and the resulting algorithms are reported in Section \ref{sec:distributed}. In that section, we also provide insights on how to establish convergence to a feasible solution or how to deduce infeasibility of the problem in a distributed manner. The convergence rate results for the proposed algorithms are described in Section \ref{sec:Convergence}. We present numerical results in Section \ref{sec:Results}, and we conclude the paper in Section \ref{sec:conclusion}.

\subsection*{Notation}
We denote the set of real $m \times n$ matrices by $\mathbb R^{m \times n}$, and $\mathbb{N}_p$ denotes the set of positive integers
$\{1,2,\ldots,p\}$. Given a set $J \subset \{1,2,\ldots,n\}$, the
matrix $E_J \in \mathbb{R}^{|J|\times n}$ is the $0$-$1$ matrix that
is obtained from an identity matrix of order $n$ by deleting the rows
indexed by $\mathbb{N}_n \setminus J$. Also, $|J|$ denotes the number of elements in set $J$. This means that $E_Jx$ is a vector with the components of $x$ that correspond to the elements in
$J$, and we denote this vector with $x_J$. Given a vector $x$, we denote the $i$th component of this vector with $x_i$. The distance from a point $x \in \mathbb{R}^n$ to a set $C
\subseteq \mathbb{R}^n$ is denoted as $\dist(x,C)$, and it is defined as
\begin{equation}
\dist(x,C) = \infimum_{y \in C} \left\| x-y \right\|.
\end{equation}
where $\| . \|$ denotes the 2-norm. Similarly, the distance between two sets $C_1, \ C_2 \subseteq
\mathbb{R}^n$  is defined as
\begin{equation}
\dist(C_1,C_2) = \infimum_{y \in C_1,  x \in C_2} \left\| x-y \right\|.
\end{equation}
The relative interior of a set $C$ is denoted $\relint(C)$, and $D =
\diag(a_1,\ldots,a_n)$ is a diagonal matrix of order $n$ with diagonal
entries $D_{ii} = a_i$. Given vectors $x^k$ for $k= 1, \dots, N$, the column vector $(x^1, \dots, x^N)$ is all of the given vectors stacked. We finally denote the so-called effective domain of a convex function, $f$, with $\dom f = \{ x \ | \ f(x) < \infty \}$.

\section{Decomposition and convex feasibility}\label{sec:Decomposition}

Given $N$ closed convex sets $\mathcal C_1,\ldots, \mathcal C_N$, a general convex feasibility problem is defined as
\begin{subequations}\label{eq:SCFP}
\begin{align}
& \find \hspace{14mm} v\\
& \subject \quad v \in \mathcal C_i, \ \ i = 1, \dots, N,
\end{align}
\end{subequations}
where $v \in \mathbb{R}^n$. We are particularly interested in the case where the description of each constraint set $\mathcal C_i$ is only dependent on a small subset of the variables in the vector $v$. Let us denote the ordered set of indices of variables that appear in the description of the $i$th constraint by $J_i$. We also denote  the ordered set of indices of constraints for which their description depend on $v_i$ by $\mathcal I_i$, i.e., $\mathcal{I}_i = \{k \ | \ i \in J_k\}$. We call a CFP loosely coupled if $|\mathcal I_i|\ll N$ for all $i = 1, \dots, n$. There are different ways of formulating the problem in~\eqref{eq:SCFP}, which allow us to design several new or use various already existing algorithms for solving this problem. In order to unify the analysis of such algorithms, we utilize so-called error bounds, which are the subject of the next subsection.
\subsection{Error Bounds and Bounded Linear Regularity}\label{sec:Error}

Error bounds quantify the distance to the solution set of a problem and they become zero only when we arrive at a solution of the problem. The use of error bounds is common in analysis of iterative algorithms, \cite{pan:97}. For CFPs, authors in \cite{bet:03} and \cite{gub:99} consider the use of
\begin{align}
T(v) = \maximum_i \left\{\dist(v,\mathcal C_i)\right\},
\end{align}
as the error bound, when analyzing projection-based algorithms. Note that $T(v) = 0$ if and only if $v \in \bigcap_{i=1}^{N}\mathcal C_i$. Based on this error bound, the closed convex sets, $\mathcal C_i$, for $i= 1, \dots, N$, are said to be boundedly linearly regular, if  for every bounded set $B$ there exists $\theta_B>0$ such that
\begin{align}\label{eq:LEB}
\forall \ v \in B \quad \dist\left(v,\bigcap_{i=1}^{N} \mathcal C_i\right) \leq \theta_B \maximum_i \left\{\dist\left(v,\mathcal C_i\right)\right\}.
\end{align}
This allows us to bound the distance to the intersection of these sets, which is very difficult or expensive to compute, by $T(v)$ which can usually be calculated easily, \cite{bet:03}. It was shown by Bauschke~\etal\ \cite{bau:99} and Beck \& Teboulle \cite{bet:03} that Slater's condition for a CFP implies bounded linear regularity, i.e., for a general CFP where $\mathcal C_1, \ldots, \mathcal C_k$ are polyhedral sets and $\mathcal C_{k+1},\ldots, \mathcal C_{N}$ are general closed, convex sets, \eqref{eq:LEB} holds if
\begin{align}
  \left( \bigcap_{i=1}^k \mathcal C_i \right) \bigcap \left( \bigcap_{i=k+1}^N
    \relint(\mathcal C_i) \right) \neq \emptyset.
 \label{e-slater-cfp}
\end{align}
Bounded linear regularity proves to be essential in the analysis of the proposed algorithms in this paper. Next, we investigate some of the approaches for solving the feasibility problem in~\eqref{eq:SCFP}.

\subsection{Projection Algorithms and Convex Feasibility Problems}\label{sec:Projection}

One of the possible approaches for solving the CFP in \eqref{eq:SCFP}, is to neglect the coupling structure among the constraint sets, and use projection algorithms for finding a solution in the intersection of $N$ convex sets. Among such projection algorithms, cyclic projection algorithm (CPA), maximum distance projection algorithm (MDPA) and mean projection algorithm (MPA) are some of the most widely used ones, where only MPA is suitable for solving \eqref{eq:SCFP} in a distributed manner. This follows from the fact that at each iteration MPA uses
\begin{align}\label{eq:cimmino}
  v^{(k+1)} := \sum_{i=1}^N \alpha_i^{(k)} P_{\mathcal C_i}(v^{(k)})
\end{align}
for updating the iterates, where $\sum_{i=1}^N \alpha_i^{(k)} = 1$ and $\alpha_1^{(k)},\ldots,\alpha_N^{(k)} > 0$. Notice that the updating procedure in \eqref{eq:cimmino} is highly parallelizable. That is because the projections can be performed in parallel and simultaneously. Assuming $N$ computing agents, each agent $i$ then computes $P_{\mathcal C_i}(v^{(k)})$, and communicates with all the other agents to update the iterate as in \eqref{eq:cimmino}. Hence this algorithm, requires global communication among all the agents. In \cite{bet:03}, it was shown that in case the sets in \eqref{eq:SCFP} are boundedly linearly regular, the algorithm enjoys a linear rate of convergence, where
\begin{align}\label{eq:cimcon}
  \dist\left(v^{(k+1)}, \bigcap_{i=1}^{N} \mathcal C_i\right) \leq \gamma_B \dist\left(v^{(k)}, \bigcap_{i=1}^{N} \mathcal C_i\right),
\end{align}
with
\begin{align}
  \gamma_B = \sqrt{1 - \frac{\minimum_i \ \{ \alpha_i^{(k)} \}}{\theta_B^2}},
\end{align}
where $\theta_B>0$ depends on the starting point $x^{(0)}$. This dependence follows from the fact that $\theta_B$ should satisfy \eqref{eq:LEB} with $B = \{v \,|\, \|v - z\| \leq \|v^{(0)} - z\|\}$ for any $z \in \mathcal C$. Iusem \& De~Pierro, \cite{iue:86}, have proposed an accelerated variant of this algorithm that takes as the next iterate a convex combination of the projections of $v^{(k)}$ on only the sets for which $v^{(k)} \notin \mathcal C_i$. This generally improves the rate of convergence when only a few constraints are violated. However, neglecting the structure in \eqref{eq:SCFP} can drastically deteriorate the performance of this algorithm, \cite{kho:13}. Also in case Slater's condition is not satisfied, e.g., when the problem is infeasible,~\eqref{eq:cimcon} does not hold and this algorithm can perform arbitrarily bad, \cite{bet:03}. In the upcoming subsections, we show how the structure in the coupling among the constraints in \eqref{eq:SCFP} can be exploited, which allows us to reformulate the problem in other ways.

\subsection{Decomposition and Product Space Formulation}\label{sec:Product}

Having the structure in the constraints in \eqref{eq:SCFP} in mind, we define $N$ lower-dimensional sets
\begin{align}
  \bar{\mathcal C}_i = \{ s^i \in \mathbb{R}^{|J_i|} \,|\, E_{J_i}^T s^i \in\mathcal C_i
  \}, \quad i=1,\ldots, N,
\end{align}
such that $s^i \in \bar{\mathcal C}_i$ implies $E_{J_i}^T s^i \in\mathcal C_i$. This allows us to rewrite the standard form CFP in \eqref{eq:SCFP} as
\begin{subequations}\label{e-cfp-decomposed}
  \begin{align}
   & \mbox{find}  \hspace{14mm}      \ s^1,s^2,\ldots,s^N,v  \\
    &\subject  \hspace{2mm} \ s^i \in \bar{\mathcal C}_i, \quad i=1,\ldots, N \\
                      & \hspace{20mm} \ s^i = E_{J_i} v, \quad i=1, \ldots, N\label{eq:consensus}
  \end{align}
\end{subequations}
where the equality constraints are the so-called coupling or global consensus constraints that ensure that the local variables $s^1,\ldots,s^N$ are consistent with one another. In other words, if the constraints $v \in \mathcal C_i$ and $v \in \mathcal C_j$ ($i\neq j$) both involve $v_k$, then the $k$th component of $E_{J_i}^T s^i$ and $E_{J_j}^T s^j$ must be equal. This formulation decomposes the so-called global variable $v$ into $N$ coupled local variables $s^1,\ldots, s^N$. This allows us to rewrite the problem as a CFP with two sets
\begin{align}\label{eq:PCFP}
  \begin{array}{ll}
    \mbox{find}     & \ \  S   \\
    \subject & \ \ S \in \mathcal C, \  S \in \mathcal D
  \end{array}
\end{align}
where
\begin{align*}
 S &= (s^1,\ldots,s^l) \in \mathbb{R}^{|J_1|}\times\cdots \times\mathbb{R}^{|J_l|} \\
\mathcal C & = \bar{\mathcal C}_1\times \cdots \times \bar{\mathcal C}_l \\
\mathcal D & = \{ \bar{E}v\,|\, v \in \mathbb{R}^n\}  \\
\bar{E} & = \begin{bmatrix} E_{J_1}^T   &\cdots &
  E_{J_l}^T \end{bmatrix}^T.
\end{align*}
The formulation \eqref{eq:PCFP} can be thought of as a ``compressed'' product space formulation of a CFP as described in~\eqref{eq:SCFP}, and it is similar to the consensus optimization problems described in \cite[Sec. 7.2]{boyd:11}, \cite[Sec. 3.4]{ber:97}. The problem in~\eqref{eq:PCFP} can now be solved using  von~Neumann's and Dykstra's alternating projections (AP) methods, which are methods for finding solutions in the intersection of two sets.

\subsubsection{Von Neumann's alternating projections}\label{sec:Von}
Given the two sets, $\mathcal C$ and $\mathcal D$, and a starting point $v^{(0)}$, von~Neumann's AP method computes two sequences
\begin{subequations}\label{eq:vonNeumann}
\begin{align}
  S^{(k+1)} &= P_{\mathcal{C}}\left( V^{(k)} \right)\\
  V^{(k+1)} &= P_{\mathcal{D}}\left( S^{(k+1)} \right).
\end{align}
\end{subequations}
where $V^{(k)} = \bar E v^{(k)}$. If the CFP in \eqref{eq:PCFP} is feasible, i.e., $\mathcal C \cap \mathcal D \neq \emptyset$, then both sequences converge to a point in $\mathcal C \cap \mathcal D$, \cite{bet:03,bau:93}. The updates in \eqref{eq:vonNeumann}, result in the following iterative algorithm
 \begin{subequations}\label{eq:vonNeumannP}
\begin{align}
\label{eq:vonNeumannP1}
S^{(k+1)} &=P_{\mathcal{C}}\left( V^{(k)} \right)\\
\notag
&= \left( P_{\bar{\mathcal C}_1}\left( E_{J_1}v^{(k)} \right), \ldots, P_{\bar{\mathcal C}_N}\left( E_{J_N}v^{(k)} \right) \right)\\
\label{eq:vonNeumannP2}
V^{(k+1)} &= \bar{E}\underbrace{\left( \bar{E}^T\bar{E} \right)^{-1}\bar{E}^T S^{(k+1)}}_{v^{(k+1)}},
\end{align}
\end{subequations}
where \eqref{eq:vonNeumannP1} and \eqref{eq:vonNeumannP2} are projections onto $\mathcal C$ and onto the column space of $\bar E$, respectively. Note that the projection onto the set $\mathcal C$ can be computed in parallel by $N$ computing agents, i.e., agent $i$ computes $s_i^{(k)} =
P_{\bar{\mathcal C}_i}(E_{J_i}v^{(k)})$, and the second projection can be interpreted as a consensus step that can be solved via distributed averaging. The details of a distributed implementation of \eqref{eq:vonNeumannP} are later discussed in Section \ref{sec:distributed}.
In case the sets $\mathcal C$ and $\mathcal D$ are boundedly linearly regular, it follows from \cite[Cor. 2.1]{bet:03} that
\begin{align}\label{eq:voncon}
  \dist\left(S^{(k+1)}, \bigcap_{i=1}^{N} \mathcal C_i\right) \leq \gamma_B \dist\left(S^{(k)}, \bigcap_{i=1}^{N} \mathcal C_i\right)
\end{align}
with
\begin{align}
  \gamma_B = \sqrt{1 - \frac{1}{\theta_B^2}}
\end{align}
where $\theta_B > 0$ depends on the starting point as is the case for~\eqref{eq:cimcon} of MPA. It was shown in \cite{bau:93,bau:94}, that in case the problem in \eqref{eq:PCFP} is infeasible
\begin{align}\label{eq:Infeasibility}
  V^{(k)} - S^{(k)}, \ V^{(k)} - S^{(k+1)} \rightarrow d, \quad \|d \| = \dist(\mathcal{C},\mathcal{D}),
\end{align}
where since $\mathcal C$ is assumed to be closed, $\dist(\mathcal{C},\mathcal{D})$ is attained. Theoretically, this result provides the possibility to detect infeasibility of \eqref{eq:PCFP} by monitoring the sequences in~\eqref{eq:Infeasibility}. However, to the best knowledge of the authors, the rates of convergence of the sequences $ V^{(k)} - S^{(k)}$ and $V^{(k)} - S^{(k+1)}$ to $d$ or $\|d \|$ to $\dist(\mathcal{C},\mathcal{D})$ have not yet been established.

\subsubsection{Dykstra's alternating projections}\label{sec:Dykstra}
The CFP in \eqref{eq:PCFP} can also be solved using Dykstra's AP method, where
\begin{subequations}\label{eq:Dykstra}
\begin{align}
S^{(k+1)} &= P_{\mathcal{C}}(V^{(k)}-\bar\lambda^{(k)})\\
\label{eq:Dykstra1}
V^{(k+1)} &= P_{\mathcal{D}}(S^{(k)}) \\
\bar\lambda^{(k+1)}&=\bar \lambda^{(k)}+(S^{(k+1)}-V^{(k+1)}).
\end{align}
\end{subequations}
and $\bar \lambda = (\bar\lambda^1, \dots, \bar\lambda^N)$. Note that this algorithm is a special case of Dykstra's AP method where one of the sets is affine, \cite{bau:94}. Similar to von~Neumann's AP method, in case $\mathcal C \cap \mathcal D \neq \emptyset$ the iterates $V^{(k)}$ and $S^{(k)}$ converge to a point in $\mathcal C \cap \mathcal D$, \cite{bau:94}. The updates in~\eqref{eq:Dykstra}, result in the following iterative algorithm
\begin{subequations}\label{eq:DykstraT}
\begin{align}
\label{eq:DykstraT1}
S^{(k+1)} &= P_{\mathcal{C}}(V^{(k)}-\bar\lambda^{(k)})\\
\notag
&= \left( P_{\bar{C}_1}\left( E_{J_1}v^{(k)} - \bar \lambda^{1,(k)} \right), \ldots, P_{\bar{C}_N}\left( E_{J_N}v^{(k)}- \bar \lambda^{N,(k)} \right) \right)\\
\label{eq:DykstraT2}
V^{(k+1)} &= P_{\mathcal{D}}(S^{(k)}) \\
\bar \lambda^{(k+1)}&= \bar \lambda^{(k)}+(S^{(k+1)}-V^{(k+1)}).
\end{align}
\end{subequations}
As can be seen from \eqref{eq:DykstraT1}, this algorithm is also highly parallelizable. This is discussed in more detail in Section~\ref{sec:distributed}. Unlike von~Neumann's AP method, the iterative algorithm in~\eqref{eq:Dykstra} does not necessarily converge with a linear rate even when the underlying sets are boundedly linearly regular. Similar to von Neumann's AP method, in case the CFP in \eqref{eq:PCFP} is infeasible the sequences $ V^{(k)} - S^{(k)}$ and $V^{(k)} - S^{(k+1)}$ converge to $d$, however, their rates of convergence are not known, \cite{bau:94}.

\subsection{Convex Minimization Formulation}\label{sec:Convex}

The problem in \eqref{eq:PCFP}, can also be reformulated as a convex minimization problem. Let $\mathcal I_{\mathcal C}(S)$ and $\mathcal I_{\mathcal D}(S)$ be the indicator functions for the sets $\mathcal C$ and $\mathcal D$, where an indicator function for a set, e.g., $\mathcal A$, is defined as
\begin{equation}\label{eq:Indicator}
\mathcal I_{\mathcal A}(x) =
\begin{cases}
\infty \hspace{15mm} x \not\in \mathcal A\\
0 \hspace{17mm} x \in \mathcal A
\end{cases}.
\end{equation}
and hence $\dom (\mathcal I_{\mathcal A}) = \mathcal A$. The CFP in \eqref{eq:PCFP} can then be equivalently rewritten as the following convex minimization problem
\begin{align}\label{eq:indmin}
\minimize_{S} \quad \mathcal I_{\mathcal C}(S) + \mathcal I_{\mathcal D}(S).
\end{align}
Despite the equivalence between the problems in \eqref{eq:PCFP} and~\eqref{eq:indmin}, the minimization problem is not defined in case the CFP in~\eqref{eq:PCFP} is infeasible, since the effective domain of the cost function would then be empty. This limits our capability to draw conclusions regarding the infeasibility of the corresponding CFP using this formulation. In order to circumvent this issue, we define the following unconstrained minimization problems

\begin{align}\label{eq:ConvexReform1}
    \minimize_{S}  \quad F_1(S) := \frac{1}{2}\sum_{i=1}^{N} \|Ês^i - P_{\bar{\mathcal C}_i}(s^i) \|^2  + \mathcal I_{\mathcal D}(S),
\end{align}
and
\begin{align}\label{eq:ConvexReform2}
    \minimize_{S} \quad F_2(S) := \frac{1}{2}\sum_{i=1}^{N} \|Ês^i - P_{\bar{\mathcal C}_i}(s^i) \|^2  +  \frac{1}{2}\ \|ÊS - P_{\mathcal D}(S) \|^2,
\end{align}
which are well-defined even when the CFP in \eqref{eq:PCFP} is infeasible. Note that the problems in \eqref{eq:ConvexReform1} and \eqref{eq:ConvexReform2} are not equivalent to the CFP in \eqref{eq:PCFP}. However, these minimization problems always have at least one solution and admit an optimal solution with zero objective value if and only if the problem in \eqref{eq:PCFP} is feasible. In fact the optimal solution then constitutes a solution for \eqref{eq:PCFP}. Similarly, the minimization problems in \eqref{eq:ConvexReform1} and~\eqref{eq:ConvexReform2}, yield a non-zero optimal objective value if and only if the CFP in~\eqref{eq:PCFP} is infeasible. In the upcoming sections, we explain how these minimization problems facilitate the design of distributed algorithms for solving the CFP in~\eqref{eq:PCFP}.


\section{Proximity Operators and Proximal Splitting}\label{sec:prox}
Consider the problem of minimizing a sum of $p$ closed convex functions
\begin{align}\label{eq:prox}
\minimize \quad F(x) =  f_1(x) + \dots + f_p(x).
\end{align}
Through the use of their so-called proximity operators, \cite{com:11}, proximal splitting algorithms allow us to perform this minimization by considering each of the terms in the cost function separately. Proximity operators are defined as follows.
\begin{definition}[\cite{com:11}]
Given a closed convex function $f :Ê \mathbb{R}^n \rightarrow \mathbb{R}$, then for every $x \in \mathbb{R}^n$, the proximity operator of the function $f$, $\text{prox}_f : \mathbb{R}^n \rightarrow \mathbb{R}^n $, is defined as the unique minimizer of the following optimization problem,
\begin{align*}
\minimize_y \quad f(y) + \frac{1}{2} \| x-y \|^2.
\end{align*}
\end{definition}

Keeping in mind the problems in \eqref{eq:ConvexReform1} and  \eqref{eq:ConvexReform2}, we only consider the case where the cost function consists of two terms, i.e., $p = 2$. Depending on the characteristics of the terms in the cost function, we are allowed to employ different proximal splitting algorithms. Next, we review some of the most widely used of such methods.
\subsection{Forward-backward Splitting}\label{sec:fb}
The forward-backward proximal splitting algorithm is suitable for cases where at least one of the two terms in the cost function is differentiable with a Lipschitz continuous gradient. Assume that $f_1$ and $f_2$ are both closed convex functions and that the problem
\begin{align}\label{eq:2p}
\minimize_x \quad f_1(x) + f_2(x)
\end{align}
has at least one solution. Let $f_1(x)$ be differentiable. In this case, the problem can be solved using Algorithm \ref{alg:alg1}.
\begin{algorithm}[H]
\caption{Forward-backward method \cite{com:11,com:05}}\label{alg:alg1}
\begin{algorithmic}[1]
\small
\State{Given $\varepsilon \in (0,\minimum(1,1/L)]$ and $x^{(1)}$}
\For {$k  = 1, 2 \dots$}
\State{$\gamma^{(k)} \in [ \varepsilon , 2/L - \varepsilon ]$}
\State{$\lambda^{(k)} \in [ \varepsilon, 1]$}
\State $y^{(k+1)} = x^{(k)} - \gamma^{(k)} \nabla f_1(x^{(k)})$
\State $x^{(k+1)} = (1-\lambda^{(k)})x^{(k)} + \lambda^{(k)} \text{prox}_{\gamma^{(k)}f_2}(y^{(k+1)})$
\EndFor
\normalsize
\end{algorithmic}
\end{algorithm}
In this algorithm $\gamma^{(k)}$ is the so-called gradient step, $\lambda^{(k)}$ is a relaxation parameter and $L$ is the Lipschitz constant of $\nabla f_1$. It was shown in \cite{lev:66}, \cite[Thm. 3.1]{bec:09}, that if $\gamma^{(k)} < 2/L$ and $\lambda^{(k)}= 1$,
\begin{align}\label{eq:fbcon}
F(x^{(k)}) - F(x^{\ast}) \leq \frac{L\| x^{(0)} - x^{(k)} \|^2}{2k},
\end{align}
where $x^\ast$ is any optimal solution for the problem in \eqref{eq:prox}. It is also possible to obtain better rate of convergence of the objective value by combining this algorithm with the so-called 1-memory accelerated gradient methods, \cite{com:11}. This comes at the expense of a more complicated algorithm. Let $l(x;y) = {f_1}(y) + \left<  \nabla f_1(y),x-y \right> + f_2(x)$ and define
\begin{align*}
D(x,y) = h(x) - h(y) - \left<  \nabla h(y),x-y \right> ,
\end{align*}
where $h$ is a strictly convex function. The general format for 1-memory accelerated gradient methods can then be presented as in Algorithm \ref{alg:alg2}, \cite{tse:08}.
\begin{algorithm}[H]
\caption{1-memory accelerated gradient method \cite{tse:08}}\label{alg:alg2}
\begin{algorithmic}[1]
\small
\State{Given $\theta^{(1)} \in (0,1]$ and $x^{(1)}, g^{(1)}$}
\For {$k  = 1, 2 \dots$}
\State $y^{(k+1)} = (1-\theta^{(k)})x^{(k)} + \theta^{(k)}g^{(k)}$
\State {$g^{(k+1)} = \argmin_x \left\{  l\left(x;y^{(k+1)}\right) + \theta^{(k)}LD\left(x,g^{(k)}\right) \right\}$}
\State $\hat x^{(k+1)} = (1-\theta^{(k)})x^{(k)} + \theta^{(k)}g^{(k+1)}$
\State{Choose $x^{(k+1)}$ to be no worse than $\hat x^{(k+1)}$ in $l(x;y^{(k+1)}) + L/2\|x - y^{(k+1)} \|^2$}
\State{Choose $\frac{1-\theta^{(k+1)}}{\left(\theta^{(k+1)}\right)^2} \leq \frac{1}{\left(\theta^{(k)}\right)^2}$}
\EndFor
\normalsize
\end{algorithmic}
\end{algorithm}
Depending on the choice of function $h(\cdot)$, and how we choose to compute $x^{(k)}$ and $\theta^{(k)}$, we end up in different accelerated gradient methods, \cite{bec:11,tse:08}. In case we choose $h(x) = \frac{1}{2} \| x - y\|^2$ and merge the fifth and sixth steps of Algorithm \ref{alg:alg2} by choosing
\begin{align*}
x^{(k+1)} = \argmin_x \{  l(x;Y^{(k+1)}) + \frac{L}{2}\| x- Y^{(k+1)} \|^2 \},
\end{align*}
we can summarize the combination of the forward-backward splitting algorithm with 1-memory accelerated gradient method as Algorithm \ref{alg:alg3}.
\begin{algorithm}[H]
\caption{Accelerated forward-backward method}\label{alg:alg3}
\begin{algorithmic}[1]
\small
\State{Given $\theta_1 \in (0,1]$ and $x^{(1)}, g^{(1)}$}
\For {$k  = 1, 2 \dots$}
\State $y^{(k+1)} = (1-\theta^{(k)})x^{(k)} + \theta^{(k)}g^{(k)}$
\State $b^{(k+1)} = g^{(k)}-\frac{1}{\theta^{(k)}L}\nabla f_1(y^{(k+1)})$
\State {$g^{(k+1)} = \text{prox}_{\frac{1}{\theta^{(k)}L}f_2}\left( b^{(k+1)} \right)$}
\State $c^{(k+1)} = y^{(k+1)}-\frac{1}{L}\nabla f_1(y^{(k+1)})$
\State{$x^{(k+1)} = \text{prox}_{\frac{1}{L}f_2}\left(c^{(k+1)} \right)$}
\State{Choose $\frac{1-\theta^{(k+1)}}{\left(\theta^{(k+1)}\right)^2} \leq \frac{1}{\left(\theta^{(k)}\right)^2}$}
\EndFor
\normalsize
\end{algorithmic}
\end{algorithm}
There are different convergence results for such algorithms which are dependent on different choices of $\theta^{(k)}$, e.g., see \cite[Cor. 1]{tse:08} and \cite[Thm. 4.4]{bec:09}, where all suggest rates of convergence of order $\mathcal O(1/k^2)$ of the objective value function. In other words
\begin{align}\label{eq:afbcon}
F(x^{(k)}) - F(x^\ast) \leq \mathcal O (\frac{1}{k^2}).
\end{align}
%
\subsection{Splitting Using Alternating Linearization Methods}\label{sec:ADMM}
The splitting of the problem in \eqref{eq:prox} for $p=2$, can also be performed by introducing auxiliary constraints as
\begin{equation}\label{eq:complex}
\begin{split}
&\minimize_{x,y} \quad f_1(x) + f_2(y)\\
&\subject  \hspace{3mm} x=y.
\end{split}
\end{equation}
Assume that $f_1$ and $f_2$ are both differentiable with Lipschitz continuous gradients. Such linear equality constrained optimization problems can then be solved using the so-called alternating linearization method (ALM), \cite{gol:12}. Let
\begin{align*}
Q^{f_2}_{\mu}(x,y) &= f_1(x) + f_2(y) + \left< \nabla f_2(y),x-y\right> + \frac{1}{2\mu}\| x-y \|^2\notag\\
& = f_1(x) + f_2(y) + \frac{1}{2\mu}\left\| x-(y-\mu\nabla f_2(y)) \right\|^2,
\end{align*}
\begin{align*}
Q^{f_1}_{\mu}(y,x) &= f_2(y) + f_1(x) + \left< \nabla f_1(x),y-x\right> + \frac{1}{2\mu}\| x-y \|^2\notag\\
& = f_2(y) + f_1(x) + \frac{1}{2\mu}\left\| x-(y+\mu\nabla f_1(x)) \right\|^2,
\end{align*}
and define the following update rules
\begin{equation}\label{eq:ALM_1}
\begin{split}
x^{(k+1)} & = \argmin_{x} \quad Q^{f_2}_{\mu_1}(x,y^{(k)}) \\
 &= \text{prox}_{\mu_1 f_1}(y^{(k)}-\mu\nabla f_2(y^{(k)})),
\end{split}
\end{equation}
and
\begin{equation}\label{eq:ALM_2}
\begin{split}
y^{(k+1)} & = \argmin_{y} \quad Q^{f_1}_{\mu_2}(y,x^{(k+1)}) \\
 &= \text{prox}_{\mu_2 f_2}(x^{(k+1)}-\mu\nabla f_1(x^{(k+1)})).
\end{split}
\end{equation}
The ALM scheme for solving the problem in~\eqref{eq:complex} can then be written as in Algorithm~\ref{alg:alg6}.
\begin{algorithm}[H]
\caption{ALM}\label{alg:alg6}
\begin{algorithmic}[1]
\small
\State{Given $\mu_1, \mu_2 > 0$ and $y^{(1)}$}
\For {$k  = 1, 2 \dots$}
\State $x^{(k+1)} =  \text{prox}_{\mu_1 f_1}(y^{(k)}-\mu\nabla f_2(y^{(k)}))$
\State $y^{(k+1)} = \text{prox}_{\mu_2 f_2}(x^{(k+1)}-\mu\nabla f_1(x^{(k+1)}))$
\EndFor
\normalsize
\end{algorithmic}
\end{algorithm}
Goldfarb \etal \  in \cite[Cor. 2.4]{gol:12} showed that in case $0<\mu_1\leq 1/L_1$ and $0<\mu_2\leq 1/L_2$, where $L_1$ and $L_2$ are Lipschitz constants for the gradients of $f_1$ and $f_2$ respectively,
\begin{align}\label{eq:ALMcon}
F(y^{(k)}) - F(x^{\ast}) \leq \frac{\| x^{(1)} - x^{\ast} \|^2}{2(\mu_1+\mu_2)k}, \quad \forall \ k>1.
\end{align}
In the same paper, the authors also propose an accelerated variant of ALM, which is reported in Algorithm \ref{alg:Falg6}.
\begin{algorithm}[H]
\caption{Fast ALM \cite{gol:12}}\label{alg:Falg6}
\begin{algorithmic}[1]
\small
\State{Given $\mu_1, \mu_2 > 0$, $t^{(1)} = 1$ and $z^{(1)}=y^{(1)}$}
\For {$k  = 1, 2 \dots$}
\State $x^{(k+1)} =  \text{prox}_{\mu_1 f_1}(z^{(k)}-\mu_1\nabla f_2(z^{(k)}))$
\State $y^{(k+1)} = \text{prox}_{\mu_2 f_2}(x^{(k+1)}-\mu_2\nabla f_1(x^{(k+1)}))$
\State $ t^{(k+1)} = \frac{1 + \sqrt{1+t^{(k)2}}}{2}$
\State $ z^{(k+1)} = y^{(k+1)} + \frac{t^{(k)}-1}{t^{(k+1)}}( y^{(k+1)} - y^{(k)})$
\EndFor
\normalsize
\end{algorithmic}
\end{algorithm}
It was shown in \cite[Cor. 3.5]{gol:12}, that in case $\mu_1$ and $\mu_2$ are chosen in the same manner as for the ALM algorithm,
\begin{align}\label{eq:FALMcon}
F(y^{(k)}) - F(x^{\ast}) \leq \frac{2\| x^{(1)} - x^{\ast} \|^2}{(\mu_1+\mu_2)k^2}, \quad \forall \ k>1
\end{align}
Notice that Algorithm \ref{alg:Falg6} is similar to Algorithm~\ref{alg:alg3}. This can be particularly observed by comparing the steps 5, 7, 8 and~1 in Algorithm~\ref{alg:alg3} with steps  3, 4, 5 and 6 in Algorithm~\ref{alg:Falg6}, respectively.
\begin{rem}
The alternating linearization method is very similar to the so-called alternating direction method of multipliers (ADMM), \cite{ber:97,boyd:11}, and in fact Algorithm \ref{alg:Falg6} is equivalent to a symmetric variant of ADMM, \cite{gol:12}. We have chosen not to investigate ADMM or its other variants, since most of their convergence rate results rely on strong convexity of at least one of the terms in the cost function, \cite{den:12,gold:12}, which is not the case for neither of the problems in \eqref{eq:ConvexReform1} and \eqref{eq:ConvexReform2}.
\end{rem}

\subsection{Douglas-Rachford Splitting}\label{sec:Douglas}
In case $\ridom f_1 \ \cap \ \ridom f_2 \neq \emptyset$ and if the problem in \eqref{eq:2p} has at least one solution, we can use the so-called Douglas-Rachford splitting algorithm for solving the optimization problem. Note that unlike forward-backward splitting, this method does not require any of the objective function terms to be differentiable. The scheme for solving the optimization problem using this method is described in Algorithm~\ref{alg:alg4}.
\begin{algorithm}[H]
\caption{Douglas-Rachford method \cite{com:11,com:07}}\label{alg:alg4}
\begin{algorithmic}[1]
\small
\State{Given $\varepsilon \in (0,1), \gamma>0$ and $y^{(1)}$}
\For {$k  = 1, 2 \dots$}
\State{$\lambda^{(k)} \in [ \varepsilon, 2-\varepsilon]$}
\State $x^{(k+1)} = \text{prox}_{\gamma f_1}(y^{(k)})$
\State $y^{(k+1)} = y^{(k)} + \lambda^{(k)}\left( \text{prox}_{\gamma f_2}(2x^{(k+1)}-y^{(k)})-x^{(k+1)}\right)$
\EndFor
\normalsize
\end{algorithmic}
\end{algorithm}
%
\begin{rem}
Douglas-Rachford splitting is one of the principle classes of splitting methods and includes many splitting methods as special cases. Particularly, it was shown in \cite{gol:89,gab:83} that ADMM and ALM fall within this class of splitting methods. The convergence of this splitting method (and its variants) has been studied thoroughly in the literature, \cite{eck:92,law:87,com:07}. However, due to its generality, the convergence rate of this splitting method in its most general format is yet to be established. Hence, although we utilize this algorithm for solving CFPs, this limits our capability to provide convergence rate results for the resulting algorithm.
\end{rem}
Next, we will apply the proximal splitting algorithms described in this section to the problems in \eqref{eq:ConvexReform1} and \eqref{eq:ConvexReform2}.
\section{Distributed Solution}\label{sec:distributed}
In this section, we propose several distributed algorithms that can be used to solve the feasibility problem in \eqref{eq:SCFP}. In sections~\ref{sec:PDI} and \ref{sec:APD}, we describe and discuss the distributed algorithms, and in Section~\ref{sec:feas}, we investigate how the convergence of these methods can be established in a distributed manner when the problem in~\eqref{eq:PCFP} is either feasible or infeasible.

\subsection{Proximal splitting and distributed implementation}\label{sec:PDI}

In order to facilitate providing a distributed solution for the feasibility problem in \eqref{eq:SCFP} and due to the similarity between the problems in \eqref{eq:ConvexReform1}, \eqref{eq:ConvexReform2} and~\eqref{eq:prox}, we employ proximal splitting algorithms. We apply algorithms \ref{alg:alg1} and \ref{alg:alg3} for solving the minimization problem in \eqref{eq:ConvexReform1}. To be able to use these algorithms, we identify $f_1(S)$ as $\frac{1}{2}\sum_{i=1}^{N}\|Ês^i - P_{\bar{\mathcal{C}}_i}(s^i) \|^2$ and $f_2(S)$ as $\mathcal I_{\mathcal D}(h)$. The proximity operators for these functions are  given as
\begin{subequations}
\begin{align}
\text{prox}_{f_1}(S) &= \frac{S+P_{\mathcal C}(S)}{2} \\
\text{prox}_{f_2}(S) &= P_{\mathcal D}(S),
\end{align}
\end{subequations}
Note that the proximity operator computation of $f_1(S)$ is highly parallelizable. Assume that a network of $N$ computing agents is available. Then $\text{prox}_{f_1}$ can be computed in a distributed manner where each of the $N$ agents calculates $\left(s^i+P_{\bar{\mathcal C}_i}(s^i)\right)/2$ individually. Considering the definition of the set $\mathcal D$, the proximity operator of $f_2$ is merely a linear projection and is given as
\begin{align}
\text{prox}_{f_2}(S) = \bar E (\bar E^T\bar E)^{-1}\bar E^T S.
\end{align}
Note that
\begin{align*}
\bar{E}^T \bar{E} = \diag( |\mathcal{I}_1|, \ldots, |\mathcal{I}_{N}|),
\end{align*}
and hence, $\text{prox}_{f_2}$ describes the required communication and interaction between the agents in the network. For instance, define $b=(\bar E^T\bar E)^{-1} \bar E^T$, then each of the components of $b$  can be expressed as
\begin{align}\label{eq:ComponentUpdate}
b_j = \frac{1}{|\mathcal{I}_j|}\sum_{q \in \mathcal{I}_j}^{} \left( E_{J_q}^Ts^{q,(k+1)} \right)_j.
\end{align}
 As a result in order to compute this quantity, the agents in the set $\mathcal I_{j}$ must interact with one another. In other words, this requires each agent $i$ to communicate with all the agents in
\begin{align}
\text{Ne}(i) = \left\{ j \ | \ J_i\cap J_j \neq \emptyset  \right\},
\end{align}
which are referred to as neighbours of agent $i$. This interpretation of $\text{prox}_{f_2}$ is later used in the description of the proposed algorithms for distributed feasibility analysis. In order to solve the minimization problem in \eqref{eq:ConvexReform2} we use algorithms \ref{alg:alg6}--\ref{alg:alg4}, where in this case, $f_2(S)$ is identified as $\frac{1}{2}\|ÊS - P_{\mathcal D}(S) \|^2$. The proximity operator for this function is given as
\begin{align}\label{eq:case2}
\text{prox}_{f_2}(S) &= \frac{S+P_{\mathcal D}(S)}{2},
\end{align}
The proximity operator of $f_2$ for this case, is also dependent on computing projections onto the consensus set. Hence, similar to the previous case, $\text{prox}_{f_2}$ will also describe the communication and interaction among agents.
\subsubsection{Forward-backward algorithm}\label{sec:Projfb}

Considering the description of the functions $f_1$ and $f_2$ in problem \eqref{eq:ConvexReform1}, $f_1$ is differentiable and
\begin{align*}
\nabla f_1(S) =  S- P_{\mathcal C}(S),
\end{align*}
which is Lipschitz continuous with Lipschitz constant $L=1$.
Applying the forward-backward method to this problem results in the following update rules
\begin{align*}
Y^{(k+1)} &= S^{(k)} -\gamma^{(k)}\left( S^{(k)} - P_{\mathcal C}(S^{(k)}) \right)\\
& = (1- \gamma^{(k)})S^{(k)} + \gamma^{(k)} P_{\mathcal C}(S^{(k)})
\end{align*}
\begin{align*}
S^{(k+1)} &= (1-\lambda^{(k)})S^{(k)} + \lambda^{(k)} P_{\mathcal D}(Y^{(k+1)})\\ &=(1-\lambda^{(k)})S^{(k)} +\lambda^{(k)} \bar E (\bar E^T\bar E)^{-1}\bar E^T Y^{(k+1)},
\end{align*}
\begin{algorithm}
\caption{Forward-backward method}\label{alg:alg7}
\begin{algorithmic}[1]
\small
\State{Given $\varepsilon \in (0,1]$, $v^{(1)}$ and $S^{(1)} = \bar E v^{(1)}$}
\For {$k  = 1, 2 \dots$}
\State{$\gamma^{(k)} \in [ \varepsilon , 2 - \varepsilon ]$}
\State{$\lambda^{(k)} \in [ \varepsilon, 1]$}
\For {$i  = 0, 1 \dots, N$}
\State{$y^{i,(k+1)} = (1- \gamma^{(k)})s^{i,(k)} + \gamma^{(k)} P_{\bar{\mathcal C}_i}(s^{i,(k)})$}
\State {Communicate with all agents $r$ belonging to $\text{Ne}(i)$}
\For {all $j \in J_i$}
 \State $v_j^{(k+1)} = \frac{1}{|\mathcal{I}_j|}\sum_{q \in \mathcal{I}_j}^{} \left( E_{J_q}^T y^{q,(k+1)} \right)_j$
  \EndFor
\State{$s^{i,(k+1)} = (1-\lambda^{(k)})s^{i,(k)} +   \lambda^{(k)}v^{(k+1)}_{J_i}$}
 \EndFor
 \EndFor
 \normalsize
 \end{algorithmic}
 \end{algorithm}
where $Y^{(k+1)} = (y^{1,(k+1)}, \dots, y^{N,(k+1)})$. Note that, if we choose $S^{(1)} = \bar E v^{(1)}$ then $S^{(k)} \in \mathcal D$ for all $k>1$. The resulting distributed solution based on the forward-backward algorithm is summarized in Algorithm \ref{alg:alg7}. Notice that, for a constant $\lambda = 1$, this algorithm is equivalent to the two point projection method in \cite{bet:03}. Furthermore if $\gamma^{(k)}= 1$, this algorithm is von~Neumann's AP method.
\subsubsection{Accelerated forward-backward algorithm}\label{sec:afb}
It is possible to obtain faster convergence rates, by employing the accelerated forward-backward splitting to solve the problem in~\eqref{eq:ConvexReform1}. Let $Y^{(k)} = (y^{1,(k)}, \dots, y^{N,(k)})$ and $G^{(k)}$,  $U^{(k)}$, $Z^{(k)}$ be defined similarly. By applying this algorithm to the problem in \eqref{eq:ConvexReform1}, we arrive at the following update rules
\begin{subequations}\
\begin{align}
Y^{(k+1)} &= (1-\theta^{(k)}) S^{(k)}+\theta^{(k)}G^{(k)}\label{eq:up1}\\
U^{(k+1)} &= G^{(k)} - \frac{1}{\theta^{(k)}}\left( Y^{(k+1)} - P_{\mathcal C}(Y^{(k+1)}) \right)\label{eq:up2}\\
G^{(k+1)} & = P_{\mathcal D}( U^{(k+1)}) \label{eq:up3}\\
Z^{(k+1)} & = P_{\mathcal C}(Y^{(k+1)})\\
S^{(k+1)} & = P_{\mathcal D}( Z^{(k+1)})
\end{align}
\end{subequations}
Note that, similar to the forward-backward algorithm, in case we choose $S^{(1)} = G^{(1)} = \bar E v^{(1)}$, we will have $Y^{(k)}, G^{(k)}, S^{(k)} \in \mathcal D$ for all $k \geq 1$. Substituting \eqref{eq:up2} into~\eqref{eq:up3} and by using \eqref{eq:up1}, we can then simplify these update rules as follows
\begin{align*}
Y^{(k+1)} &= (1-\theta^{(k)}) S^{(k)}+\theta^{(k)}G^{(k)}\\
G^{(k+1)} & = \bar E (\bar E^T\bar E)^{-1}\bar E^T U^{(k+1)}\\
& = \bar E (\bar E^T\bar E)^{-1}\bar E^T \left( G^{(k)} - \frac{1}{\theta^{(k)}}\left( Y^{(k+1)} - P_{\mathcal C}(Y^{(k+1)}) \right)  \right)\\
& = G^{(k)} - \frac{1}{\theta^{(k)}} \bar E (\bar E^T\bar E)^{-1}\bar E^T \\ &\hspace{20mm}\left( (1-\theta^{(k)}) S^{(k)}+\theta^{(k)}G^{(k)} - P_{\mathcal C}(Y^{(k+1)}) \right)\\
& = \frac{\theta^{(k)} - 1}{\theta^{(k)}}S^{(k)} + \frac{1}{\theta^{(k)}} \bar E (\bar E^T\bar E)^{-1}\bar E^TP_{\mathcal C}(Y^{(k+1)})\\
S^{(k+1)} & = \bar E (\bar E^T\bar E)^{-1}\bar E^T P_{\mathcal C}(Y^{(k+1)})
\end{align*}
The resulting distributed algorithm can then be summarized as  in Algorithm~\ref{alg:alg8}.
\begin{algorithm}[H]
\caption{Accelerated proximal gradient method}\label{alg:alg8}
\begin{algorithmic}[1]
\small
\State{Given $\theta^{(0)} \in (0,1]$, $v^{(0)}$ and $G^{(1)}=S^{(1)} = \bar E v^{(0)}$}
\For {$k  = 1, 2, \dots$}
\For {$i  = 1, 2 \dots, N$}
\State{$y^{i,(k+1)} = (1- \theta^{(k)})s^{i,(k)} + \theta^{(k)}g^{i,(k)}$}
\State {Communicate with all agents $r$ belonging to $\text{Ne}(i)$}
\For {all $j \in J_i$}
 \State $v_j^{(k+1)} = \frac{1}{|\mathcal{I}_j|}\sum_{q \in \mathcal{I}_j}^{} \left( E_{J_q}^T P_{\bar{\mathcal C}_q}(y^{q,(k+1)}) \right)_j$
 \EndFor
 \State{$g^{i,(k+1)} = \frac{\theta^{(k)} - 1}{\theta^{(k)}}s^{i,(k)} + \frac{1}{\theta^{(k)}}v_{J_i}^{(k+1)}$}
 \State{$s^{i,(k+1)} = v_{J_i}^{(k+1)}$}	
 \EndFor
 \State{Choose $\theta^{(k+1)}$ such that $\frac{1-\theta^{(k+1)}}{\theta^{(k+1)2}}\leq \frac{1}{\theta^{(k)2}}$}
 \EndFor
 \normalsize
 \end{algorithmic}
\end{algorithm}

\subsubsection{ALM}\label{sec:ProjALM}

It is also possible to devise a distributed feasibility analysis algorithm by applying ALM to the formulation in~\eqref{eq:ConvexReform2}. Define $\nu^{(k)} = (\nu^{1,(k)}, \dots, \nu^{N,(k)}) = \nabla f_1(S^{(k)})$ and $\xi^{(k)} = (\xi^{1,(k)}, \dots, \xi^{N,(k)}) = \nabla f_2(Y^{(k)})$. From the optimality conditions for \eqref{eq:ALM_1} and \eqref{eq:ALM_2}, we have
\begin{align*}
&\nabla f_1(S^{(k+1)}) + \frac{1}{\mu_1} (S^{(k+1)} - Y^{(k)}) + \nabla f_2(Y^{(k)}) = 0\\
&\nabla f_2(Y^{(k+1)}) - \frac{1}{\mu_2} (S^{(k+1)} - Y^{(k+1)}) + \nabla f_1(S^{(k+1)}) = 0
\end{align*}
which results in the following update rules for $\nu^{(k)}$ and $\xi^{(k+1)}$
\begin{equation}\label{eq:ALMup}
\begin{split}
\nu^{(k+1)} &= -\xi^{(k)} - \frac{1}{\mu_1}(S^{(k+1)}-Y^{(k)})\\
\xi^{(k+1)} &= -\nu^{(k+1)} + \frac{1}{\mu_2}(S^{(k+1)}-Y^{(k+1)})
\end{split}
\end{equation}
Applying Algorithm \ref{alg:alg6} to the problem in \eqref{eq:ConvexReform2} then results in
\begin{subequations}
\begin{align}
S^{(k+1)} &= \text{prox}_{\mu_1 f_1}(Y^{(k)}-\mu_1\xi^{(k)})\notag\\
& = \frac{1}{\mu_1 +1} \left( Y^{(k)}-\mu_1\xi^{(k)} + \mu_1 P_{\mathcal C}(Y^{(k)}-\mu_1\xi^{(k)}) \right)\\
\nu^{(k+1)} &= -\xi^{(k)} - \frac{1}{\mu_1}(S^{(k+1)}-Y^{(k)})\label{eq:ALM1}\\
Y^{(k+1)} &= \text{prox}_{\mu_2 f_2}(S^{(k+1)}-\mu_2\nu^{(k+1)})\notag\\
& = \frac{1}{\mu_2 +1} \left( S^{(k+1)}-\mu_2\nu^{(k+1)} + \mu_2 P_{\mathcal D}(S^{(k+1)}-\mu_2\nu^{(k+1)}) \right)\notag\\
& = \frac{1}{\mu_2 +1} \left( S^{(k+1)}-\mu_2\nu^{(k+1)} + \mu_2 \bar E (\bar E^T\bar E)^{-1}\bar E^T(S^{(k+1)}-\mu_2\nu^{(k+1)}) \right)\\
\xi^{(k+1)} &= -\nu^{(k+1)} + \frac{1}{\mu_2}(S^{(k+1)}-Y^{(k+1)}),\label{eq:ALM2}
\end{align}
\end{subequations}
which is obtained by combining the update rules in \eqref{eq:ALM_1}, \eqref{eq:ALM_2} and \eqref{eq:ALMup}. Since the Lipschitz constants for both $f_1$ and $f_2$ are equal to~1, we can also choose $\mu_1 = \mu_2 = 1$. The resulting distributed feasibility algorithm can then be summarized in Algorithm \ref{alg:alg10}.
\begin{algorithm}[H]
\caption{Alternating linearization method}\label{alg:alg10}
\begin{algorithmic}[1]
\footnotesize
\State{Given $Y^{(1)}$ and $\xi^{(1)} = Y^{(1)} - P_{\mathcal D}(Y^{(1)})$}
\For {$k  = 1, 2, \dots$}
\For {$i  = 1, 2 \dots, N$}
\State{$s^{i,(k+1)} = \frac{1}{2} \left( y^{i,(k)} - \xi^{i,(k)}+ P_{\bar{\mathcal C}_i}(y^{i,(k)} - \xi^{i,(k)})\right)$}
\State{$\nu^{i,(k+1)} = -\xi^{i,(k)} - (s^{i,(k+1)}-y^{i,(k)})$}
\State {Communicate with all agents $r$ belonging to $\text{Ne}(i)$}
\For {all $j \in J_i$}
 \State $v_j^{(k+1)} = \frac{1}{|\mathcal{I}_j|}\sum_{q \in \mathcal{I}_j}^{} \left( E_{J_q}^T (s^{q,(k+1)}-\nu^{q,(k+1)}) \right)_j$
 \EndFor
\State{$y^{i,(k+1)} = \frac{1}{2} \left( s^{i,(k+1)}-\nu^{i,(k+1)} +  v_{J_ i}^{(k+1)} \right)$}
\State{$\xi^{i,(k+1)} = -\nu^{i,(k+1)} + (s^{i,(k+1)}-y^{i,(k+1)})$}
 \EndFor
 \EndFor
 \normalsize
 \end{algorithmic}
\end{algorithm}
\begin{rem}
The authors in \cite{gol:12}, propose another variant of Algorithm \ref{alg:alg6} that allows for non-smooth terms in the cost function, with similar convergence results. This enables us to use ALM for solving the formulation in \eqref{eq:ConvexReform2} of the CFP. However, doing so recovers von Neumann's AP method.
\end{rem}
\subsubsection{Fast ALM}\label{sec:ProjFALM}
Following the ideas from the derivation of Algorithm \ref{alg:alg10}, we can also apply fast ALM to the problem in~\eqref{eq:ConvexReform2} as follows. Define $\nu^{(k)}= \nabla f_1(S^{(k)})$, $\xi^{(k)} = \nabla f_2(Y^{(k)})$ and similarly $\beta^{(k)} = \nabla f_2(Z^{(k)})$. From the optimality conditions of the $3^{\text{rd}}$ and $4^{\text{th}}$ steps of Algorithm \ref{alg:Falg6}, we have
\begin{align*}
&\nabla f_1(S^{(k+1)}) + \frac{1}{\mu_1} (S^{(k+1)} - Z^{(k)}) + \nabla f_2(Z^{(k)}) = 0\\
&\nabla f_2(Y^{(k+1)}) - \frac{1}{\mu_2} (S^{(k+1)} - Y^{(k+1)}) + \nabla f_1(S^{(k+1)}) = 0.
\end{align*}
Also recall that $\nabla f_2(Z^{(k)}) = Z^{(k)} - P_{\mathcal D}(Z^{(k)}) = Z^{(k)} - \bar E (\bar E^T\bar E)^{-1}\bar E^T Z^{(k)} $ is linear with respect to the input argument, and hence, by the $6$th step of Algorithm \ref{alg:Falg6}, we arrive at the following update rules
\begin{equation}\label{eq:FALMup}
\begin{split}
\nu^{(k+1)} &= -\beta^{(k)} - \frac{1}{\mu_1}(S^{(k+1)}-Z^{(k)})\\
\xi^{(k+1)} &= -\nu^{(k+1)} + \frac{1}{\mu_2}(S^{(k+1)}-Y^{(k+1)})\\
\beta^{(k+1)} & = \xi^{(k+1)} + \frac{t^{(k)}-1}{t^{(k+1)}}( \xi^{(k+1)} - \xi^{(k)})
\end{split}
\end{equation}
Combining these with the resulting update rules obtained from applying Algorithm \ref{alg:Falg6} to the problem in \eqref{eq:ConvexReform2}, yields
\begin{subequations}
\begin{align}
S^{(k+1)} &= \text{prox}_{\mu_1 f_1}(Z^{(k)}-\mu_1\beta^{(k)})\notag\\
& = \frac{1}{\mu_1 +1} \left( Z^{(k)}-\mu_1\beta^{(k)} + \mu_1 P_{\mathcal C}(Z^{(k)}-\mu_1\beta^{(k)}) \right)\\
\nu^{(k+1)} &= -\beta^{(k)} - \frac{1}{\mu_1}(S^{(k+1)}-Z^{(k)})\label{eq:ALM1}\\
Y^{(k+1)} &= \text{prox}_{\mu_2 f_2}(S^{(k+1)}-\mu_2\nu^{(k+1)})\notag\\
& = \frac{1}{\mu_2 +1} \left( S^{(k+1)}-\mu_2\nu^{(k+1)} + \mu_2 P_{\mathcal D}(S^{(k+1)}-\mu_2\nu^{(k+1)}) \right)\notag\\
& = \frac{1}{\mu_2 +1} \left( S^{(k+1)}-\mu_2\nu^{(k+1)} + \mu_2 \bar E (\bar E^T\bar E)^{-1}\bar E^T(S^{(k+1)}-\mu_2\nu^{(k+1)}) \right)\\
\xi^{(k+1)} &= -\nu^{(k+1)} + \frac{1}{\mu_2}(S^{(k+1)}-Y^{(k+1)}),\label{eq:ALM2}\\
t^{(k+1)} &= \frac{1 + \sqrt{1+t^{(k)2}}}{2}\\
Z^{(k+1)} &= Y^{(k+1)} + \frac{t^{(k)}-1}{t^{(k+1)}}( Y^{(k+1)} - Y^{(k)})\\
\beta^{(k+1)} & = \xi^{(k+1)} + \frac{t^{(k)}-1}{t^{(k+1)}}( \xi^{(k+1)} - \xi^{(k)}),
\end{align}
\end{subequations}
where similar to the previous algorithm, we can choose $\mu_1 = \mu_2 = 1$. This algorithm is summarized in Algorithm \ref{alg:alg11}.
\begin{algorithm}[H]
\caption{Fast alternating linearization method}\label{alg:alg11}
\begin{algorithmic}[1]
\footnotesize
\State{Given $Z^{(1)} = Y^{(1)}$, $\beta^{(1)} = Z^{(1)} - P_{\mathcal D}(Z^{(1)})$ and $t^{(1)} = 1$}
\For {$k  = 1, 2, \dots$}
\For {$i  = 1, 2 \dots, N$}
\State{$s^{i,(k+1)} = \frac{1}{2} \left( z^{i,(k)} - \beta^{i,(k)}+  P_{\bar{\mathcal C}_i}(z^{i,(k)} - \beta^{i,(k)})\right)$}
\State{$\nu^{i,(k+1)} = -\beta^{i,(k)} - (s^{i,(k+1)}-z^{i,(k)})$}
\State {Communicate with all agents $r$ belonging to $\text{Ne}(i)$}
\For {all $j \in J_i$}
 \State $v_j^{(k+1)} = \frac{1}{|\mathcal{I}_j|}\sum_{q \in \mathcal{I}_j}^{} \left( E_{J_q}^T (s^{q,(k+1)}-\nu^{q,(k+1)}) \right)_j$
 \EndFor
\State{$y^{i,(k+1)} = \frac{1}{2} \left( s^{i,(k+1)}-\nu^{i,(k+1)} + v_{J_ i}^{(k+1)} \right)$}
\State{$\xi^{i,(k+1)} = -\nu^{i,(k+1)} +(s^{i,(k+1)}-y^{i,(k+1)})$}
\State $ t^{(k+1)} = \frac{1 + \sqrt{1+t^{(k)2}}}{2}$
\State $ z^{i,(k+1)} = y^{i,(k+1)} + \frac{t^{(k)}-1}{t^{(k+1)}}( y^{i,(k+1)} - y^{i,(k)})$
\State $\beta^{i,(k+1)}  = \xi^{i,(k+1)} + \frac{t^{(k)}-1}{t^{(k+1)}}( \xi^{i,(k+1)} - \xi^{i,(k)})$
 \EndFor
 \EndFor
 \normalsize
 \end{algorithmic}
\end{algorithm}
\begin{rem}
In \cite{gol:12}, another variant of Algorithm \ref{alg:Falg6} is suggested that can handle non-differentiable terms in the cost function and can deliver similar convergence rate results. This variant of the algorithm includes a skipping step which does not allow an efficient distributed implementation and requires global communication of iterates among all agents.
\end{rem}
\subsubsection{Douglas-Rachford algorithm}\label{sec:ProjDouglas}

We now apply the Douglas-Rachford algorithm to the minimization problem in~\eqref{eq:ConvexReform2}. This results in the following update rules,

\begin{align*}
S^{(k+1)}& = \frac{1}{\gamma+1}\left(Y^{(k)} + \gamma P_{\mathcal C}(Y^{(k)}) \right)\\
Y^{(k+1)} & = Y^{(k)} + \lambda^{(k)}\left( \text{prox}_{\gamma f_2} (2S^{(k+1)}-Y^{(k)}) - S^{(k+1)}\right)\\
& = Y^{(k)}  +  \lambda^{(k)}\\ & \hspace{12mm}\left( \frac{1}{\gamma +1}\left( 2S^{(k+1)}-Y^{(k)} + \gamma P_{\mathcal D}(2S^{(k+1)}-Y^{(k)}) \right) - S^{(k+1)}\right)\\
& = Y^{(k)}  +  \lambda^{(k)}\\ & \hspace{8mm}\left( \frac{2}{\gamma +1}S^{(k+1)}-\frac{1}{\gamma +1}Y^{(k)} + \frac{\gamma}{\gamma +1} P_{\mathcal D}(2S^{(k+1)}-Y^{(k)}) - S^{(k+1)}\right)\\
& = Y^{(k)}  +  \lambda^{(k)}\left( \frac{1-\gamma}{\gamma +1}S^{(k+1)}-\frac{1}{\gamma +1}Y^{(k)} + \frac{\gamma}{\gamma +1} P_{\mathcal D}(2S^{(k+1)}-Y^{(k)})\right)
\end{align*}
 The resulting iterative algorithm is reported in Algorithm~\ref{alg:alg9}.
\begin{algorithm}[H]
\caption{Douglas-Rachford method}\label{alg:alg9}
\begin{algorithmic}[1]
\footnotesize
\State{Given $\varepsilon \in (0,1)$, $\gamma > 0$ and $Y^{(1)}$}
\For {$k  = 1, 2, \dots$}
\State{$\lambda^{(k)} \in [ \varepsilon , 2 - \varepsilon ]$}
\For {$i  = 1, 2 \dots, N$}
\State{$s^{i,(k+1)} = \frac{1}{\gamma + 1} \left( y^{i,(k)}+ \gamma P_{\bar{\mathcal C}_i}(y^{i,(k)})\right)$}
\State{$z^{i,(k+1)} = 2s^{i,(k+1)} - y^{i,(k)} $}
\State {Communicate with all agents $r$ belonging to $\text{Ne}(i)$}
\For {all $j \in J_i$}
 \State $v_j^{(k+1)} = \frac{1}{|\mathcal{I}_j|}\sum_{q \in \mathcal{I}_j}^{} \left( E_{J_q}^T z^{q,(k+1)} \right)_j$
 \EndFor
\State{$y^{i,(k+1)} = y^{i,(k)} + \lambda^{(k)}\left( \frac{1-\gamma}{\gamma +1}s^{i,(k+1)} -\frac{1}{\gamma +1}y^{i,(k)} + \frac{\gamma}{\gamma +1}v_{J_i}^{(k+1)} \right)$}
 \EndFor
 \EndFor
 \normalsize
 \end{algorithmic}
\end{algorithm}
Note that this algorithm is  similar to the method proposed method in \cite{par:11} for large-scale distributed learning.
\subsection{Distributed Implementation of von~Neumann's and Dykstra's AP method}\label{sec:APD}

Similar to the algorithms presented in sections~\ref{sec:Projfb}--\ref{sec:ProjDouglas}, it is also possible to implement von~Neumann's and Dykstra's AP methods in a distributed manner. The distributed version of these algorithms are presented in algorithms \ref{alg:alg12} and~\ref{alg:alg13}.
\begin{algorithm}[H]
\caption{Von~Neumann's AP method}\label{alg:alg12}
\begin{algorithmic}[1]
\small
\State{Given $x^{(1)}$}
\For{k = 1, 2, \dots}
\For{i = 1, 2, \dots,N}
\State $s^{i,(k+1)} =  P_{\bar{\mathcal C}_i}\left( v_{J_i}^{(k)}\right).$
\State {Communicate with all agents $r$ belonging to $\text{Ne}(i)$}
\For {all $j \in J_i$}
 \State $v_j^{(k+1)} = \frac{1}{|\mathcal{I}_j|}\sum_{q \in \mathcal{I}_j}^{} \left( E_{J_q}^T s^{q,(k+1)} \right)_j$
 \EndFor
\EndFor
\EndFor
\normalsize
 \end{algorithmic}
 \end{algorithm}
\begin{algorithm}[H]
\caption{Dykstra's AP method}\label{alg:alg13}
\begin{algorithmic}[1]
\small
\State{Given $x^{(1)}$ and $\bar{\lambda}^{(1)} = 0$}
\For{k = 1, 2, \dots}
\For{i = 1, 2, \dots,N}
\State $s^{i,(k+1)} =  P_{\bar{\mathcal C}_i}\left(  v_{J_i}^{(k)} -\bar{\lambda}^{i,(k)} \right)$
\State {Communicate with all agents $r$ belonging to $\text{Ne}(i)$}
\For {all $j \in J_i$}
 \State $v_j^{(k+1)} = \frac{1}{|\mathcal{I}_j|}\sum_{q \in \mathcal{I}_j}^{} \left( E_{J_q}^T s^{q,(k+1)} \right)_j$
 \EndFor
 \State $\bar{\lambda}^{i,(k+1)} = \bar{\lambda}_i^{(k)} + \left( s^{i,(k+1)}-v_{J_i}^{(k+1)} \right)$
\EndFor
\EndFor
\normalsize
 \end{algorithmic}
 \end{algorithm}

\subsection{Local convergence tests}\label{sec:feas}

In case \eqref{eq:PCFP} is feasible, algorithms \ref{alg:alg7}--\ref{alg:alg13} converge to a feasible solution of the problem. Local convergence tests check the convergence of the iterates to a feasible solution or detect infeasibility of the problem in a distributed manner with minimal communication. Unlike the so-called global tests, which require transmission of the local variables to a central unit, local methods only demand each agent to merely declare its local variables feasibility or convergence status with respect to its local constraints and/or objective value. Recall that for feasible problems applying the proposed algorithms to the formulations in \eqref{eq:ConvexReform1} and \eqref{eq:ConvexReform2} will generate sequences that converge to an optimal solution, \cite{com:05,tse:08,com:07,bau:11,com:11}, which yield zero objective value. Also for the case~\eqref{eq:PCFP} is infeasible, all the proposed algorithms converge to a solution of the problems in \eqref{eq:ConvexReform1} and \eqref{eq:ConvexReform2}, however, this solution does not result in zero objective value. Based on the aforementioned observation, we propose an approach for establishing convergence to a feasible solution or infeasibility of the problem. This approach is based on monitoring the convergence of the objective function value and the satisfaction of local constraints.
\subsubsection{Convergence of the objective value}
One of the ways to establish convergence of the proposed algorithms is through monitoring the so-called relative change of objective value (they intend to minimize), at each iteration. In case this quantity falls below a certain threshold, we can deduce convergence of the algorithm to a solution. Particularly, given a sequence $\{ d^{(k)} \}$, the relative change of this sequence at iteration $k+1$ can be defined as
\begin{align}
R^{(k+1)} = \frac{\left\| d^{(k+1)} - d^{(k)} \right\|}{\left\| d^{(k)} \right\|}.
\end{align}
For algorithms \ref{alg:alg7} and \ref{alg:alg8} that concern the problem in \eqref{eq:ConvexReform1}, we then monitor the following quantity

\begin{align}
R_1^{(k+1)} = \frac{\left|  \left\|S^{(k+1)} - P_{\mathcal C}\left( S^{(k+1)}\right) \right\|^2 - \left\|S^{(k)} - P_{\mathcal C}\left( S^{(k)}\right) \right\|^2  \right|}{\left\|S^{(k)} - P_{\mathcal C}\left( S^{(k)}\right) \right\|^2}.
\end{align}
This is because $S^{(k)} \in \mathcal D$, $\forall \ k \geq 1$, for both algorithms. Monitoring this quantity locally, however, requires that all agents communicate their iterates to all other agents in the network. In order to alleviate this issue, we instead consider monitoring an upper bound for $R_1^{(k+1)}$. Notice that

\begin{equation}\label{eq:RC1}
\begin{split}
R_1^{(k+1)} &\leq \frac{\sum_{i=1}^N\left|\left\| s^{i,(k+1)} - P_{\mathcal{\bar C}_i}\left( s^{i,(k+1)} \right)\right\|^2 - \left\| s^{i,(k)} - P_{\mathcal{\bar C}_i}\left( s^{i,(k)} \right)\right\|^2\right|}{\left\|S^{(k)} - P_{\mathcal C}\left( S^{(k)}\right) \right\|^2}\\
&\leq \sum_{i=1}^N\frac{\left|\left\| s^{i,(k+1)} - P_{\mathcal{\bar C}_i}\left( s^{i,(k+1)} \right)\right\|^2 - \left\| s^{i,(k)} - P_{\mathcal{\bar C}_i}\left( s^{i,(k)} \right)\right\|^2\right|}{\left\| s^{i,(k)} - P_{\mathcal{\bar C}_i}\left( s^{i,(k)} \right)\right\|^2}\\ & =:  \sum_{i=1}^N R_1^{i,(k+1)}.
\end{split}
\end{equation}
As can be seen from \eqref{eq:RC1}, this upper bound can now be monitored in a distributed manner. Then if all the local relative changes, i.e., $R_1^{i,(k+1)}$, fall below a certain threshold, we can infer convergence of the algorithm to a solution. Hence, we can deduce convergence with little communication.

For algorithms \ref{alg:alg10}--\ref{alg:alg9}, which concern the problem in \eqref{eq:ConvexReform2}, the monitored relative change takes the following form

\begin{equation}
\begin{split}
R_2^{(k+1)} = \frac{\left| F_2\left( Y^{(k+1)} \right) - F_2\left( Y^{(k)} \right)\right|}{F_2\left( Y^{(k)} \right)}
\end{split}
\end{equation}
where $F_2(\cdot)$ is defined as in \eqref{eq:ConvexReform2}. The relative change $R_2^{(k+1)}$ can also be bounded in a similar manner as \eqref{eq:RC1} as follows

\begin{equation}
\begin{split}
R_2^{(k+1)} &\leq \sum_{i=1}^N\frac{\left|\| y^{i,(k+1)} - P_{\mathcal{\bar C}_i}\left( y^{i,(k+1)} \right)\|^2 - \| y^{i,(k)} - P_{\mathcal{\bar C}_i}\left( y^{i,(k)} \right)\|^2\right|}{\| y^{i,(k)} - P_{\mathcal{\bar C}_i}\left( y^{i,(k)} \right)\|^2+\| y^{i,(k)} - E_{J_i}C^{(k)}\|^2} \\ & \hspace{10mm} +  \frac{\left|\| y^{i,(k+1)} - E_{J_i}C^{(k+1)}\|^2 - \| y^{i,(k)} - E_{J_i}C^{(k)}\|^2\right|}{\| y^{i,(k)} - P_{\mathcal{\bar C}_i}\left( y^{i,(k)} \right)\|^2+\| y^{i,(k)} - E_{J_i}C^{(k)}\|^2} \\ & =:  \sum_{i=1}^N R_2^{i,(k+1)}.
\end{split}
\end{equation}
where $C^{(k)} = (E^T E)^{-1}E^TY^{(k)}$. As a result the convergence of algorithms \ref{alg:alg10}--\ref{alg:alg9} can also be established distributedly and with little communication. However, notice that for each agent to compute its local relative change at each iteration, i.e., $R_2^{i,(k+1)}$, additional communication among agents is required. This additional communication is required for computation of $C^{(k)}$.

\subsubsection{Feasibility of local constraints}
In case the CFP is feasible, all the proposed algorithms converge to a feasible solution. We can detect arrival at a feasible solution distributedly, by checking the feasibility of local constraints. If at a certain iteration, the local iterates of all agents satisfy their corresponding local constraints and if furthermore we have global consensus over the network, i.e., \eqref{eq:consensus} is satisfied, we can infer that we have converged to a feasible solution.

For algorithms \ref{alg:alg7} and \ref{alg:alg8}, the iterate $S^{(k)}$ already satisfy the global consensus constraints. Hence, the feasibility detection at each iteration (for these algorithms) requires each agent $i$ to check whether $s^{i,(k)} \in \bar{\mathcal C_i}$. In case this is satisfied for all $i = 1, \dots, N$, we can then infer arrival at a feasible solution. For algorithms \ref{alg:alg10}--\ref{alg:alg9}, however, this test is slightly more complicated. This is because the iterate $Y^{(k)}$ does not necessarily satisfy the global consensus constraints. Consequently, when $y^{i,(k)} \in \bar{\mathcal C_i}$ for all $i = 1, \dots, N$, the agents would still need to communicate with their neighbors to check whether global consensus is reached or not. Then, in case local constraints for all agents and global consensus constraints are both satisfied, we can infer arrival at a feasible solution. 

By combining the two methods for detecting convergence of the objective value and arrival at a feasible solution, we can now describe a distributed framework for establishing convergence to a solution as follows. At each iteration, each agent should
\begin{enumerate}
\item check the feasibility of its local iterates with respect to their corresponding local constraints. If all agents are locally feasible, communicate with neighbors to check the satisfaction of the global consensus constraint (this only applies to algorithms~\ref{alg:alg10}--\ref{alg:alg9});
\item check whether the local relative change has fallen below the predefined threshold.
\end{enumerate}
Then, 
\begin{itemize}
\item if condition (1) is satisfied for all agents, the algorithm has converged to a feasible solution;
\item if condition (2) is satisfied for all agents, the algorithm has converged and in case there exists an agent with non-zero local objective value, the CFP is infeasible.    
\end{itemize}
\begin{rem}
The convergence of algorithms \ref{alg:alg12} and \ref{alg:alg13}, can also be established in a similar manner. We refer to \cite{kho:13}, for details of the corresponding convergence detection framework.
\end{rem}

\section{Convergence Rate}\label{sec:Convergence}

In this section, we investigate the convergence results for the algorithms presented in Section \ref{sec:distributed}. We are particularly interested in the possibility of unifying the convergence rate results for proximal methods with the existing results for projection methods, which were discussed in Section \ref{sec:Decomposition}. The convergence rate results for projection methods, see \eqref{eq:cimcon} and~\eqref{eq:voncon}, are based on the distance of the iterates to the feasible set and are proven under the assumption that the underlying sets are boundedly linearly regular, (or that Slater's conditions are satisfied). In order to unify these results with convergence rate results for proximal splitting methods, we study the convergence of the algorithms presented in Section \ref{sec:distributed} by investigating the feasible and infeasible cases separately.
\subsection{Feasible problem}\label{sec:feasible}
Throughout this section we assume that the CFP in \eqref{eq:PCFP} is feasible and its underlying constraint sets are boundedly linearly regular, i.e., they satisfy \eqref{eq:LEB}.

\subsubsection{Forward-backward splitting}

In this subsection we focus on algorithms \ref{alg:alg7} and \ref{alg:alg8}, which are obtained by applying  forward-backward splitting to \eqref{eq:ConvexReform1}. As was mentioned in sections \ref{sec:Projfb} and \ref{sec:afb}, in these algorithms, the iterate $S^{(k)} \in \mathcal D$ for all $k\geq 1$. Hence, $\dist(S^{(k)}, \mathcal D) = 0$ and
\begin{align}\label{eq:fbcon1}
F(S^{(k)}) &= \frac{1}{2} \left\| S^{(k)} - P_{\mathcal C}(S^{(k)}) \right\|^2 + \mathcal I_{\mathcal D}(S^{(k)}) \\&= \frac{1}{2} \left\| S^{(k)} - P_{\mathcal C}(S^{(k)}) \right\|^2.
\end{align}
Assuming bounded linear regularity of the problem in \eqref{eq:PCFP}, we then have
\small
\begin{align*}
\dist\left( S^{(k)}, \mathcal C \cap \mathcal D \right) &\leq \theta_B \maximum \left\{ \dist\left( S^{(k)}, \mathcal C \right),\dist\left( S^{(k)},\mathcal D \right) \right\} \\ & = \theta_B\dist\left( S^{(k)}, \mathcal C \right)
\end{align*}
\normalsize
Consequently and by \eqref{eq:fbcon1}, for algorithms \ref{alg:alg7} and \ref{alg:alg8} we have that
\begin{align}
\dist \left( S^{(k)}, \mathcal C \cap \mathcal D \right) \leq \mathcal O(\frac{1}{\sqrt{k}}),
\end{align}
and
\begin{align}
\dist \left(S^{(k)}, \mathcal C \cap \mathcal D \right) \leq \mathcal O(\frac{1}{k}),
\end{align}
respectively.
%

\subsubsection{ALM splitting}
Recall that algorithms \ref{alg:alg10} and \ref{alg:alg11} are obtained by applying ALM and fast ALM to the problem in~\eqref{eq:ConvexReform2}. Assuming bounded linear regularity of the CFP, we then arrive at
\small
\begin{align*}
\dist^2\left( Y^{(k)}, \mathcal C \cap \mathcal D \right) &\leq \theta_B^2  \maximum \left\{ \dist^2\left( Y^{(k)}, \mathcal C \right),\dist^2\left( Y^{(k)},\mathcal D \right) \right\} \\&\leq \left\| Y^{(k)} - P_{\mathcal C}(Y^{(k)}) \right\|^2 + \left\| Y^{(k)} - P_{\mathcal D}(Y^{(k)}) \right\|^2 \\&= 2 F(Y^{(k)}).
\end{align*}
\normalsize
Note that in case the problem in \eqref{eq:PCFP} is feasible $F(S^{\ast}) = 0$ and hence by \eqref{eq:ALMcon} and \eqref{eq:FALMcon}, we have convergence rate results
\begin{align}
\dist \left( Y^{(k)}, \mathcal C \cap \mathcal D \right) \leq \mathcal O(\frac{1}{\sqrt{k}}),
\end{align}
and
\begin{align}
\dist \left(Y^{(k)}, \mathcal C \cap \mathcal D \right) \leq \mathcal O(\frac{1}{k}),
\end{align}
 for algorithms \ref{alg:alg10} and \ref{alg:alg11}, respectively.

\subsection{Infeasible problem}
The convergence results for algorithms \ref{alg:alg7}--\ref{alg:alg11}, are not affected by the feasibility or infeasibility of the CFP in~\eqref{eq:PCFP}. Hence, we can consider the cost function of the problems in~\eqref{eq:ConvexReform1} and \eqref{eq:ConvexReform2}, as a measure for detecting infeasibility. This is similar to the measures used for detecting infeasibility for von~Neumann's and Dykstra's AP methods, where the sequences $\|ÊV^{(k)} - S^{(k)} \|$ and $\|ÊV^{(k)} - S^{(k+1)} \|$ converge to $\dist(\mathcal{C},\mathcal{G})$. However, recall that the rate of convergence of the mentioned sequences is not established.
This is in contrast to algorithms  \ref{alg:alg7}--\ref{alg:alg11}, where we can use the existing results on convergence rate of the objective value of algorithms \ref{alg:alg1}--\ref{alg:Falg6}, to provide convergence results on certain residuals that can assist us in detecting infeasibilty of the problem.

\subsubsection{Forward-backward splitting}

Since for algorithms \ref{alg:alg7} and~\ref{alg:alg8}, the iterate $S^{(k)} \in \mathcal D$ for all $k \geq 1$, even when $\mathcal C \cap \mathcal D = \emptyset$, we have
\begin{align*}
F(S^{(k)}) = \frac{1}{2} \left\| S^{(k)} - P_{\mathcal C}(S^{(k)}) \right\|^2.
\end{align*}
Recall that if the problem in \eqref{eq:PCFP} is infeasible, the optimal objective value for the problem in \eqref{eq:ConvexReform1} will be nonzero. Hence we can establish the infeasibility of the problem, by monitoring the residual $\left\| S^{(k)} - P_{\mathcal C}(S^{(k)}) \right\|^2$, which will converge to $\dist^2\left( S^\ast, \mathcal C \right) = \dist^2\left( \mathcal C, \mathcal D \right)$. By the convergence results in \eqref{eq:fbcon} and \eqref{eq:afbcon}, we know that  for algorithms \ref{alg:alg7} and~\ref{alg:alg8}, the rates of convergence of this residual are of $\mathcal O(1/k)$ and $\mathcal O(1/k^2)$, respectively.
\subsubsection{ALM splitting}
For algorithms \ref{alg:alg10} and \ref{alg:alg11}, we can also draw similar conclusions. Recall that in case the problem in \eqref{eq:PCFP} is infeasible, the optimal objective value for the problem in \eqref{eq:ConvexReform2} will be nonzero. Hence, we can deduce infeasibility of the problem by monitoring the convergence of the objective value of the problem in \eqref{eq:ConvexReform2}, i.e.,
\begin{align*}
\left\| Y^{(k)} - P_{\mathcal C}(Y^{(k)})  \right\|^2 + \left\| Y^{(k)} - P_{\mathcal D}(Y^{(k)})  \right\|^2,
\end{align*}
which, by \eqref{eq:ALMcon} and \eqref{eq:FALMcon}, is known to converge with $\mathcal O(1/k)$ and $\mathcal O(1/k^2)$ convergence rates, respectively.
\begin{rem}
Among algorithms \ref{alg:alg7}--\ref{alg:alg9}, the ones based on the accelerated forward-backward and fast ALM, i.e., algorithms~\ref{alg:alg8} and \ref{alg:alg11}, have better convergence properties. However, Algorithm \ref{alg:alg8} has a practical advantage over the other. This is because, in order to detect convergence and arrival at a feasible solution, Algorithm~\ref{alg:alg11} requires more communication with neighbors and also a more sophisticated approach to do so. Hence, and purely based on the discussions in sections \ref{sec:feas} and \ref{sec:Convergence}, we expect Algorithm \ref{alg:alg8} to outperform the rest of the described methods.
\end{rem}

\section{Numerical Results}\label{sec:Results}

In this section, we apply algorithms \ref{alg:alg7}--\ref{alg:alg13} to a class of flow feasibility problems, and compare their performance when solving these problems. In section \ref{sec:FP}, we first describe the considered flow feasibility problem and we then apply the algorithms to feasible and infeasible flow problems in sections~\ref{sec:FeasFP} and \ref{sec:InfeasFP}, respectively.
\subsection{Flow Feasibility Problem}\label{sec:FP}

Let $G(\mathcal V,\mathcal E)$ be a directed graph, where $\mathcal V = \{ 1,\dots,N \}$ is the set of its vertices and $\mathcal E \subseteq \mathcal V \times \mathcal V$ is the set of its edges. Two nodes $i$ and $j$ are adjacent if $(i,j) \in \mathcal E$, and the set of the $i$th node's adjacent nodes are denoted as $\adj(i)$. 
Also let the nodes $u, o \in \mathcal V$ be the so-called source and sink nodes of the graph, respectively. Assume that we inject a flow $U$ to the source node. The flow feasibility problem then corresponds to the problem of assessing whether it is possible to relay $U$ from the source node to the sink node, by assigning flows to different edges in the graph and without violating flow constraints. These constraints mainly describe how different nodes in the network are allowed to relay the input flow from the source node to the sink node, by assigning flows to their edges. Let $f^i_j$ denote the flow assigned to the edge $(i,j) \in \mathcal E$. Notice that since $f^i_j$ and $f_i^j$ correspond to the flow within the same edge, we have
\begin{align}\label{eq:FPC}
f^i_j = f_i^j, \quad \forall \ (i,j) \in \mathcal E.
\end{align}
The constraints in the flow feasibility problem can then be expressed as follows.
\begin{enumerate}
\item The flow within each edge should be nonnegative and should not exceed its maximum capacity, i.e.,
\begin{align*}
0\leq f^i_j\leq c_{ij}, \quad \forall \ (i,j) \in \mathcal E,
\end{align*}
where $c_{ij}$ denotes the maximum capacity of the $(i,j)$ edge, and naturally $c_{ij} = c_{ji}$.
\item All nodes should satisfy the conservation of flow at all time. In other words, the sum of flows entering a node should be equal to the sum of flows leaving a node, i.e., for all nodes $i \in \mathcal V \setminus \{ u, o \}$
\begin{align*}
\sum_{j\in\adj(i)\setminus \mathcal O(i)} f^i_j = \sum_{j\in\mathcal O(i)} f^i_j,
\end{align*}
where $\mathcal O(i)$ denotes the set of $i$th nodes' adjacent nodes that receive flow from this node. For the nodes $u$ and $o$, this entails
\begin{align*}
\sum_{j\in\adj(u)\setminus \mathcal O(u)} f^u_j + U= \sum_{j\in\mathcal O(u)} f^u_j,
\end{align*}
and
\begin{align*}
\sum_{j\in\adj(o)\setminus \mathcal O(o)} f^o_j = \sum_{j\in\mathcal O(o)} f^o_j + U,
\end{align*}
respectively. 
\item  The sum of flows leaving a node should not exceed an internal nodal capacity (which could be private to the node), i.e., for all nodes $i \in \mathcal V \setminus \{ o \}$,
\begin{align*}
\sum_{j\in\mathcal O(i)} f^i_j \leq n_i,
\end{align*}
where $n_i$ is the $i$th node nodal capacity, and for the node $o$, this entails
\begin{align*}
\sum_{j\in\mathcal O(o)} f^o_j + U\leq n_o.
\end{align*}
\end{enumerate}
Having described the constraints of the flow feasibility problem, for $i\in \mathcal V \setminus \{ u,o \}$, define 
\begin{equation}
\bar{\mathcal C}_i = \left\{s^i \Bigg | \begin{split} &\sum_{j\in\adj(i)\setminus \mathcal O(i)} f^i_j = \sum_{j\in\mathcal O(i)} f^i_j\\ & \hspace{4mm}\sum_{j\in\mathcal O(i)} f^i_j \leq n_i\\ &\hspace{6mm} 0\leq f^i_j\leq c_{ij} \ \ \forall j\in \adj(i) 
\end{split}\right\}.
\end{equation}
Similarly for $i \in \{u,o\}$, define
\begin{equation}
\bar{\mathcal C}_u = \left\{s^u \Bigg | \begin{split} &\sum_{j\in\adj(u)\setminus \mathcal O(u)} f^u_j + U = \sum_{j\in\mathcal O(u)} f^u_j\\ 
& \hspace{4mm}\sum_{j\in\mathcal O(u)} f^u_j \leq n_u\\ 
&\hspace{7mm} 0\leq f^u_j\leq c_{uj} \ \ \forall j\in \adj(u) 
\end{split}\right\},
\end{equation}
and
\begin{equation}
\bar{\mathcal C}_o = \left\{s^o \Bigg | \begin{split} &\sum_{j\in\adj(o)\setminus \mathcal O(o)} f^o_j = \sum_{j\in\mathcal O(o)} f^o_j + U\\
 & \hspace{4mm}\sum_{j\in\mathcal O(o)} f^o_j + U \leq n_o\\ 
 &\hspace{7mm} 0\leq f^o_j\leq c_{oj} \ \ \forall j\in \adj(o) 
\end{split}\right\},
\end{equation}
where for each $i \in \mathcal V$, $s^i$ is the vector of flows of all the edges that are connected to the $i$th node. Notice that the sets $\bar{\mathcal C}_i$, for $i=1,\dots,N$, are decoupled, and \eqref{eq:FPC} describes the coupling among these constraints. In other words, \eqref{eq:FPC} defines the global consensus constraints. With this definition of $s^i$s and $\bar{\mathcal C}_i$s, the flow feasibility problem has the same format as the problem in \eqref{eq:PCFP}, where $\mathcal D$ is given by~\eqref{eq:FPC}. Next, we apply algorithms~\ref{alg:alg7}--\ref{alg:alg13}, to this class of flow feasibility problems.  

\subsection{Feasible Flow Problem}\label{sec:FeasFP}
\begin{figure}
\begin{center}
\includegraphics[width=10.5cm]{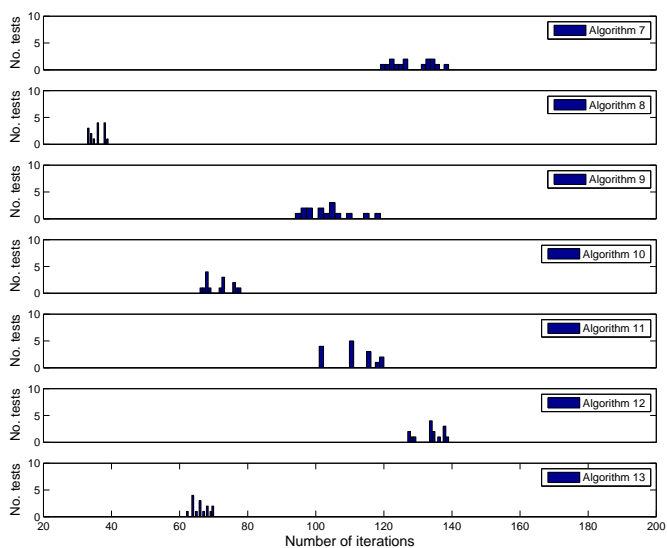}    
\caption{Number of required iterations for each algorithm to converge to a feasible solution. The figure illustrates the results achieved for all the 15 randomly generated examples.}
\label{fig:results}
\end{center}
\end{figure}      
\begin{figure}
\begin{center}
\includegraphics[width=9cm]{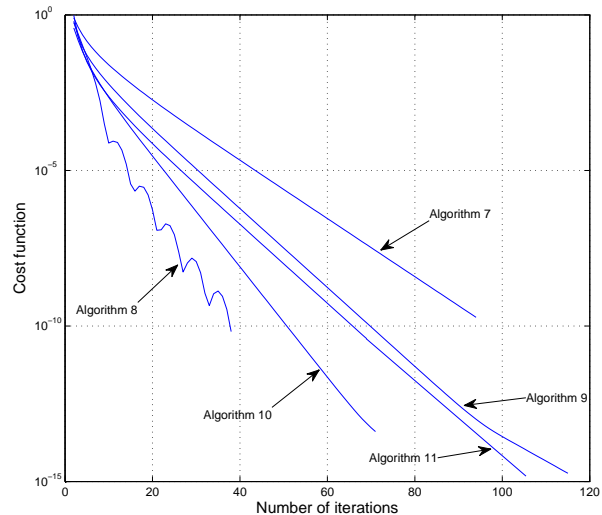}    
\caption{The evolution of the cost function being minimized for algorithms~\ref{alg:alg7}--\ref{alg:alg9}.}
\label{fig:resultsres}
\end{center}
\end{figure}
In this section, we will study and compare the performance of algorithms \ref{alg:alg7}--\ref{alg:alg13}, when they are applied to a set of flow feasibility problems. To this end, we pose flow feasibility problems over 15 connected directed graphs with 60 nodes. The source and sink nodes of each of these graphs have been chosen to be the nodes with the most number of outgoing and ingoing edges, respectively. These graphs have been generated using the algorithm presented in Appendix \ref{app:app1}. The number of variables in the feasibility problems that correspond to the generated flow problems, i.e., the number of local and global variables, varies within the range of $[936, 1146]$. Also for the sake of simplicity, we have chosen the edge capacities to be equal for all edges and also we have chosen the nodal capacities to be a proportional to the number of outgoing edges from each node, i.e., given a graph $G(\mathcal V, \mathcal E)$, $c_{ij} = \bar c$ for all $(i,j) \in \mathcal E$, and $n_i = | \mathcal O(i)| \bar n$ for all $i \in \mathcal V$. In order to assure feasibility of the generated problems, we have chosen the input flow to the source node, edge capacities and nodal capacities such that the network is capable of relaying the input flow to the sink node. Particularly, for the examples considered in this section, we have used the following procedure for finding suitable input flow, edge and nodal capacities. We first set the initial values for the input flow and the capacities as $U = 100$, $\bar c = 10$ and $n_i = | \mathcal O(i)|  \bar c/2$ for all $i=1 , \dots, N$. If the resulting flow problem is feasible, then the initial values are deemed to be suitable. Otherwise, we 
\begin{enumerate}
\item set $U := U/2$.
\item then check if the resulting flow problem is feasible; in which case we consider the current values of $\bar c$, $n_i$s and $U$ as the chosen ones. In case the problem was still infeasible, we
\begin{itemize}
\item set $\bar c := 2 \bar c$. 
\item set $n_i := | \mathcal O(i)|  \bar c/2$ for all $i=1 , \dots, N$. 
\item again check whether the resulting flow problem is feasible or not. In case the problem is feasible we have found the suitable values for $U$, $\bar c$ and $\bar n$. Otherwise, we continue by going back to step (1) of the scheme.
\end{itemize}
\end{enumerate}

For algorithms \ref{alg:alg8}--\ref{alg:alg11} that do not have any tuning parameters (except for $\theta^0$ in Algorithm \ref{alg:alg8} which is chosen to be 1) we have applied the algorithms as they are. For algorithms~\ref{alg:alg7} and \ref{alg:alg9} that do contain tuning parameters, these parameters are chosen such that the algorithms achieve their best performance for each specific example. However, we have not considered time varying parameters in these algorithms. Figure \ref{fig:results} illustrates the obtained results from this experiment. The figure shows the number of required iterations for each algorithm to converge to a feasible solution. In order to detect convergence to a feasible solution, we have utilized the proposed approach in Section~\ref{sec:feas}, where the threshold for convergence detection based on local relative changes was set to $10^{-4}$. However, for all examples, the convergence to feasible solutions was established using the condition concerning the feasibility of local constraints.  This was because this condition was satisfied prior to the convergence of the objective value. As can be seen from Figure \ref{fig:results}, Algorithm~\ref{alg:alg8} clearly outperforms the rest of the algorithms, followed by algorithms~\ref{alg:alg13} and \ref{alg:alg11}. Figure \ref{fig:resultsres}, illustrates the performance of algorithms~\ref{alg:alg7}--\ref{alg:alg9} in minimizing their corresponding cost functions, for one specific example. As was expected, algorithms~\ref{alg:alg8} and \ref{alg:alg11}, outperform the rest of the algorithms. 
\subsection{Infeasible Flow Problem}\label{sec:InfeasFP}
\begin{figure}
\begin{center}
\includegraphics[width=9cm]{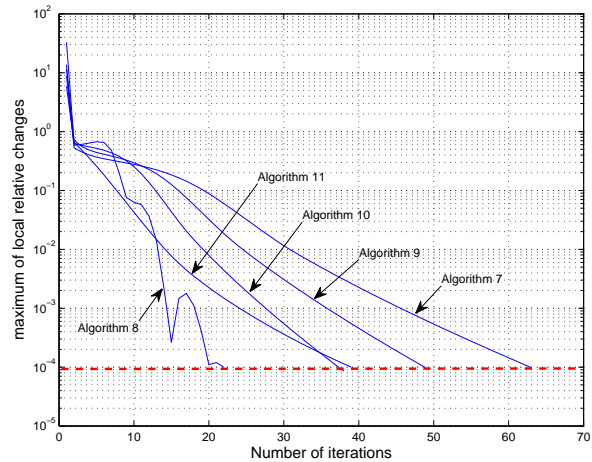}    
\caption{The maximum of local relative changes over the network for algorithms~\ref{alg:alg7}--\ref{alg:alg9} when the flow feasibility problem is infeasible.}
\label{fig:resultsIn1}
\end{center}
\end{figure}
\begin{figure}
\begin{center}
\includegraphics[width=9cm]{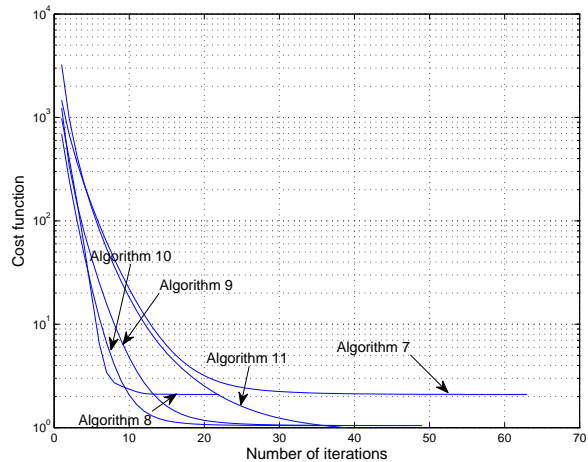}    
\caption{The evolution of the cost function being minimized for algorithms~\ref{alg:alg7}--\ref{alg:alg9} when the flow feasibility problem is infeasible.}
\label{fig:resultsIn}
\end{center}
\end{figure}
In this experiment, we consider the case when the flow feasibility problem is infeasible, where we also use a similar setup as in Section \ref{sec:FeasFP}. Particularly, we randomly generate a connected directed graph with 60 nodes, however, we design the problem such that the resulting flow problem is infeasible. To be more specific, we reduce the capacity of the edges and increase the input flow to the network, up to a point that we exceed the relaying capabilities of the network. The procedure that we used for finding suitable values for $U$, $\bar c$ and $n_i$s for this purpose, is similar to the approach discussed in Section~\ref{sec:FeasFP}. In fact, the only difference is in the steps when we change the values for $U$ and $\bar c$, where we instead set these values as $U:= 2U$ and $\bar c := \bar c/2$.

Note that, since the flow feasibility problem is infeasible, the convergence can only be established through monitoring the local relative changes. Figure
\ref{fig:resultsIn1} depicts the maximum of local relative changes across the network, for algorithms \ref{alg:alg7}--\ref{alg:alg9}. The dashed line in the figure, illustrates the threshold for convergence detection. As can be seen from the figure, Algorithm \ref{alg:alg8} outperforms the other algorithms and converges within 22 iterations. This is followed by algorithms~\ref{alg:alg11} and~Ê\ref{alg:alg9} which converged within 38  and 39 iterations, respectively. It is worth mentioning that, Dykstra's AP method also converged to a solution in 39 iterations, however, von Neumann's AP method required  892 iterations to converge. Figure \ref{fig:resultsIn} depicts the behavior of algorithms \ref{alg:alg7}--\ref{alg:alg9} when applied to an infeasible problem. Unlike the behavior that was observed in Figure~\ref{fig:resultsres}, the cost function sequence, now, converges to a nonzero constant. 
\section{Conclusions}\label{sec:conclusion}
In this paper, we presented several algorithms for solving loosely coupled convex feasibility problems distributedly and efficiently. These algorithms were the result of application of proximal splitting methods to convex minimization reformulations of product-space formulation of such CFPs. We also proposed a distributed feasibility/infeasibility detection scheme that require little communication among the agents. Furthermore, through the use of convergence rate results for proximal splitting methods, we provided a unified treatment of the convergence rate analysis of the proposed algorithms and of classical projection methods. We also studied the performance of the proposed algorithms using numerical experiments which illustrated that Algorithm \ref{alg:alg8} outperforms the rest of the algorithms. 

A possible shortcoming of the presented algorithms can be in handling infeasible problems where $\dist(\mathcal C,\mathcal D)$ is extremely small. In such cases and due to sub-linear convergence properties of these algorithms, they can require many iteration to converge (though, this was not observed in any of the 15 examples). As future research direction, we intend to further investigate this problem and devise possible remedies, e.g., by using more advanced (primal-dual) splitting methods. 
\appendices
\section{An Algorithm for Random Generation of Connected Directed Graphs}\label{app:app1}
In this appendix, we present an algorithm for generating connected and directed graphs, $G(\mathcal V, \mathcal E)$,  with $N$ vertices and with adjacency matrix  $A \in \mathbb R^{N\times N}$ given by
\begin{align}\label{eq:I}
A_{ij} = \begin{cases}
0 \quad \  \ (i,j) \notin \mathcal E \\
1 \quad \ \ (i,j) \in \mathcal E \ \text{the edge is leaving the $i$th node}\\
-1 \quad (i,j) \in \mathcal E \ \text{the edge is entering the $i$th node}
\end{cases}.
\end{align} 
Notice that, by this definition, $A$ is skew-symmetric. Next, we describe an algorithm, that allows us to randomly generate the adjacency matrix of connected and directed graphs.
\begin{algorithm}[H]
\caption{ Random Generation of Connected Directed Graphs }\label{alg:Graph}
\begin{algorithmic}[1]
\footnotesize
\State{Given $N$, $F_1 = 1$, $Iter_{m}= 1000$ and $A$ a $N \times N$ zero matrix}
\For  {$i = 1:Iter_{m}$}
\For {$i = N-1:-1:1$}
\State Set $F_2 = 1$ and $F_3 = 1$
\While  {$F_2 == 1$}
\State Generate a random $1\times i$ 0-1 vector, $x$.
\If {Number of nonzero elements in $\begin{bmatrix}A(N-i,1:N-i)&x \end{bmatrix}$ is larger than 2}
\State Set $F_2 = 0$
\EndIf
\EndWhile 
\While  {$F_3 == 1$}  
\State Randomly assign a sign to nonzero elements in $x$.
\If {There exists both positive and negative elements in $\begin{bmatrix}A(N-i,1:N-i)&x \end{bmatrix}$}
\State Set $F_3 = 0$
\State $A(N-i,N-i+1:N) = x$
\State $A(N-i+1:N,N-i) = x^T$
\EndIf
\EndWhile
\EndFor
\If {There exists both positive and negative elements in $A(N,1:N)$}
\State break
\EndIf
\EndFor
 \normalsize
 \end{algorithmic}
\end{algorithm}
Note that this algorithm can fail to generate a suitable adjacency matrix at each run, and one must continue on running the algorithm until it satisfies the conditions in the algorithm. The proposed algorithm is not an efficient method for generating connected directed graphs, and it is merely a simple methodology that we used in Section \ref{sec:Results} for generating such graphs.


\bibliographystyle{gOMS} 
\bibliography{IEEETrans}

\end{document}